\documentclass[11pt]{elsarticle}

\makeatletter
\def\ps@pprintTitle{%
 \let\@oddhead\@empty
 \let\@evenhead\@empty
 \def\@oddfoot{}%
 \let\@evenfoot\@oddfoot}
\usepackage{bm}
\usepackage{bbm}
\usepackage{amsmath}
\usepackage{amsthm}
\usepackage{empheq}
\usepackage{float}
\usepackage{lineno}
\usepackage{indentfirst}
\usepackage{cases}
\usepackage{algorithm}  
\usepackage{algorithmicx} 
\usepackage{algpseudocode} 
\usepackage{subcaption}
\usepackage{graphicx}
\usepackage{booktabs}
\usepackage{diagbox}
\usepackage[all]{xy}
\usepackage{indentfirst}
\usepackage{epstopdf}
\usepackage{epsfig}
\usepackage{overpic}
\usepackage{array}
\usepackage{color}
\usepackage{multirow}
\usepackage[mathscr]{euscript}
\usepackage{latexsym,bm}
\usepackage{tikz}
\usetikzlibrary{shapes,arrows}
\usepackage[left=2.5cm,right=2.5cm,top=2cm,bottom=2cm]{geometry}
\usepackage{enumerate}

\usepackage{ulem}

\usepackage{framed,multirow}

\usepackage{amssymb}
\usepackage{latexsym}
\usepackage{mathtools}

\usepackage{url}
\usepackage{xcolor}
\definecolor{newcolor}{rgb}{.8,.349,.1}

\renewcommand{\vec}[1]{\boldsymbol{#1}} 

\newtheorem{example}{Example}[section]

\def\bx{\bm{x}}
\def\by{\bm{y}}
\def\be{\bm{e}}
\def\me{\mathbbm{e}}
\def\D{\textup{d}}
\def\mone{\mathbbm{1}}

\makeatletter

\newcommand{\Rmnum}[1]{\expandafter@slowromancap\romannumeral #1@}
\makeatother

\tikzset{global scale/.style={
		scale=#1,
		every node/.append style={scale=#1}
	}
}
\tikzstyle{startstop} = [rectangle,rounded corners, minimum width=3cm,minimum height=1cm,text centered, draw=black,fill=red!30]
\tikzstyle{io} = [trapezium, trapezium left angle = 70,trapezium right angle=110,minimum width=3cm,minimum height=1cm,text centered,draw=black,fill=blue!30]
\tikzstyle{process} = [rectangle,minimum width=3cm,minimum height=1cm,text centered,text width =3cm,draw=black,fill=orange!30]
\tikzstyle{decision} = [diamond,minimum width=1.7cm,minimum height=0.5cm,shape aspect=2, text centered,draw=black,fill=green!30]
\tikzstyle{arrow} = [thick,->,>=stealth]

\newtheorem{assumption}{Assumption}[section]

\DeclareMathOperator*{\argmin}{\ensuremath{arg\,min}}

\begin{document}

\begin{frontmatter}

\title{\LARGE An efficient stochastic particle method for moderately high-dimensional nonlinear PDEs}

\author[PKU]{Zhengyang Lei}
\ead{leizy@stu.pku.edu.cn}

\author[PKU]{Sihong Shao}
\ead{sihong@math.pku.edu.cn}

\author[BNU]{Yunfeng Xiong}
\ead{yfxiong@bnu.edu.cn}


\address[PKU]{CAPT, LMAM and School of Mathematical Sciences,  Peking University, Beijing 100871, China}
\address[BNU]{School of Mathematical Sciences, Beijing Normal University, Beijing 100875, China}

\begin{abstract}

Numerical resolution of moderately high-dimensional nonlinear PDEs remains a huge challenge due to the curse of dimensionality for the classical numerical methods including finite difference, finite element and spectral methods. Starting from the weak formulation of the Lawson-Euler scheme, this paper proposes a stochastic particle method (SPM) by tracking the deterministic motion, random jump, resampling and reweighting of particles. 
Real-valued weighted particles are adopted by SPM to approximate the high-dimensional solution, which automatically adjusts the point distribution to intimate the relevant feature of the solution.  
A piecewise constant reconstruction with virtual uniform grid is employed to evaluate the nonlinear terms, which fully exploits the intrinsic adaptive characteristic of SPM. Combining both, SPM can achieve the goal of adaptive sampling in time. Numerical experiments on the 6-D Allen-Cahn equation and the 7-D Hamiltonian-Jacobi-Bellman equation demonstrate the potential of SPM in solving moderately high-dimensional nonlinear PDEs efficiently while maintaining an acceptable accuracy.

\end{abstract}


\end{frontmatter}

\section{Introduction}
\label{sec:intro}

Moderately high-dimensional PDEs have been widely used in modeling complicated physical and social phenomena. Typical examples include the 6-D Schr{\"o}dinger equation \cite{smyth1998numerical} or the 12-D Wigner equation \cite{XiongShao2020Overcoming} in describing the two-body quantum system, the 6-D Vlasov equation \cite{kormann2019massively} or Boltzmann equation \cite{dimarco2018efficient} for the plasma evolution, and the Hamilton–Jacobi–Bellman (HJB) equation in control theory the dimension of which is the number of players \cite{bk:Bellman1957}. In contrast to their sound theoretical validation, the notorious curse of dimensionality (CoD) poses a fundamental obstacle to numerical resolution of moderately high-dimensional problems for the cost increases exponentially with the dimensionality.

Up to now, all existing numerical attempts to alleviate CoD often fall into two distinct categories: Grid-based methods and particle-based methods. The former have made significant advances towards simulating real 6-D PDEs with the help of high-performance computing platforms, e.g., the Boltzmann equation with variable hard sphere kernel \cite{dimarco2018efficient}, the Vlasov equation \cite{kormann2019massively}, the Schr\"odinger equation for two-electron atoms \cite{smyth1998numerical}, the quantum Wigner-Coulomb system \cite{xiong2022characteristic} and the Euler-Lagrange equation  for  icosahedral quasicrystals \cite{jiang2017stability}. 
However, since the number of basis functions used still explodes with the dimension, it is entirely not trivial to extend them into PDEs in higher dimensions subject to the memory limit of current computers. When the solutions possess the low-rank structure, the tensor trains method \cite{bachmayr2023low,bachmayr2016tensor,dolgov2021tensor, tang2024solving, richter2024continuous} can be employed to address higher-dimensional problems.
Particle-based methods that utilize the radial basis functions  \cite{LiChenChen2013}, sophisticated sampling \cite{Giles2008,WangSloan2005,YanCaiZeng2013,su2023new} and optimization \cite{BayerEigelSallandtTrunschke2023} have also achieved a huge success in solving certain high-dimensional PDEs in engineering and finance, but their performance relies heavily on the  low effective dimension underlying the problems \cite{WangSloan2005}. Recently, Deep-learning-based PDE solvers, which are also recognized as particle-based methods \cite{EMaWu2019}, have emerged with a rapid development and been applied to solve high-dimensional PDEs in many areas \cite{HanJentzenE2018,RaissiPerdikarisKarniadakis2019,HurePhamWarin2020}. 
However, it is sensitive to the hyperparameter tuning \cite{kast2024positional}, and also poses some challenging problems in terms of the theoretical understanding of optimization aspects \cite{richter2021solving}.

In this paper, we propose a stochastic particle method (SPM) for solving moderately high-dimensional nonlinear PDEs, e.g., the 7-D HJB equation, which automatically distributes particles relevant to their high-dimensional solution landscapes.
The motivation of using stochastic particles to integrate PDEs dates back to the CFL paper \cite{CFL1928} where the boundary value problem of elliptic equations was connected with the random walk problem through a finite difference scheme. Whether a SPM may break CoD has attracted attentions in the recent decade. Several attempts based on the Feynman-Kac formula were proposed, such as the multi-level Monte Carlo method for linear parabolic equations \cite{Giles2008}, the branching diffusion process algorithm \cite{HenryLabordere2019} and the multi-level Picard iterative algorithm \cite{weinan2019multilevel} for semi-linear parabolic equations. However, these methods usually obtain the function value at a prescribed point every run and there still remain some key issues to be further studied. For example, the exponential growth of particles in the branching diffusion process algorithm result in a short-time simulation, and the Lipschitz constant of nonlinear terms is restricted to be small in the multi-level Picard iteration. It is more important to note that the Feynman-Kac formula expresses the solution as the expectation of a positive weighted particle system. What if particles possess different charges? This will make the mathematical tools of probability theory almost unusable \cite{mainini2012description}. Techniques from optimal transport may be helpful in studying evolution PDEs with sign-changing solutions, 
like the Wasserstein distance for signed measures \cite{ambrosio2011gradient} and the Kantorovitch-Rubinstein duality \cite{piccoli2019wasserstein}. Motivated by the CFL paper \cite{CFL1928}, which obtains stochastic particle interpretation directly from the numerical schemes rather than from the PDEs, this work proposes an alternative way to obtain SPM for moderately high-dimensional nonlinear PDEs with sign-changing solutions by utilizing the corresponding Lawson-Euler \cite{lawson1967generalized} time discretization.

When dealing with CoD, the accuracy of an efficient numerical method can be expected at current stage is only half order, like 
the Monte Carlo with convergence rate $\mathcal{O}(N^{-1/2})$ and $N$ being the sample size, rather than a higher order one.
In keeping with such spirit, we choose the Lawson-Euler scheme 
which linearizes the nonlinear terms using the solutions at the previous time step and thus only has first-order accuracy in time. After that, just like most traditional numerical methods, we propose to use the weak formulation for obtaining stochastic particle representation in high-dimensional space, the usefulness of which can be twofold. First, it allows us to conveniently extract required information from the high-dimensional unknown functions. By plugging different test functions into the weak formulation, it is able to figure out the desired information of the solutions more than the function value at a fixed point.  Second, it appears as a high-dimensional integral that depends on time and the underlying PDE establishes a bridge in the time direction. That is, SPM can be seen as using the Monte Carlo to approximate the time-dependent integral,  
and simulates the evolution of the PDE through particle motion and weight update and so on. This is the key point to automatically adjust the particles to  intimate the feature of high-dimensional solutions, thereby achieving the goal of adaptive sampling.
It should be pointed out that, the proposed SPM does not depend on the Feynman-Kac formula.

Starting from the weak formulation of the Lawson-Euler scheme,  SPM uses real-valued weighted particles to approximate the high-dimensional solution in the weak sense. It tracks the deterministic motion, random jump, resampling and reweighting of particles. To reduce the cost of reconstructing the solution values via weighted particles
and in accordance with the first-order accuracy of the Lawson-Euler scheme, it suggests to adopt a uniform piecewise constant reconstruction with virtual uniform grid (VUG) by fully taking advantage of the high adaptivity of stochastic particles. It should be noted that, finding the optimal piecewise constant reconstruction is an NP-hard problem \cite{LiYangWong2016} with high computational complexity.  SPM is fully adaptive in nature and easy to parallelize. Using the distributed parallel technology via Message Passing Interface (MPI), it solves the 6-D Allen-Cahn equation and the 7-D HJB equation efficiently while maintaining a reasonable error tolerance, and the accuracy can be systematically improved by increasing the number of particles.

SPM combines Monte Carlo with grid-based methods. To evaluate the nonlinear term, the piecewise constant reconstruction with virtual uniform grid is employed, causing the computational cost still increases exponentially with the dimension. Nevertheless, compared with the full tensor girds often used in classical numerical methods, SPM can reduce storage to a certain extent by exploiting the intrinsic adaptive characteristic that particles move with the solution (see Table~\ref{VUGtable}). Currently, SPM can only deals with nonlinear PDEs in moderately high-dimensional problems. By contrast, other attempts for high-dimensional problems including deep learning based \cite{han2017deep,nusken2021interpolating,GaoWang2023} or tensor train based methods \cite{bachmayr2023low,bachmayr2016tensor,dolgov2021tensor, tang2024solving, richter2024continuous} might potentially deal with problems in much higher dimensions. The advantage of SPM lies in that it grows on sound and solid mathematical ground, without learning and an a priori assumption that the solution has a low-rank structure.

The rest of this paper is organized as follows. In Section \ref{sec1}, we introduce the theory of stochastic particle method. Section \ref{implement section} presents more implementation details. Some 1-D benchmark tests are conducted in Section \ref{sec:acc}. Section \ref{sec:ac} and Section \ref{sec:hjb} show numerical experiments for solving the 6-D Allen-Cahn equation and the 7-D HJB equation, respectively. 
The paper is concluded in Section \ref{sec:con} with a few remarks.

\section{The stochastic particle method}
\label{sec1}

Consider nonlinear PDEs of the following form,
\begin{equation}\label{cauchy problem}
\begin{aligned}
	\frac{\partial } {\partial t}  u(\bm{x}, t) &= \mathcal{L} u(\bx, t) +  f(t,\bx,u,\vec{\nabla }u), \quad \bx \in \mathbb{R}^d, ~ t> 0, \\
	u(\bm{x}, 0) &= u_0(\bx), \quad \bx \in \mathbb{R}^d,
\end{aligned}
\end{equation}
where $\mathcal{L}$ represents the linear operator, 
$f(t,\bx,u,\vec{\nabla }u)$ gives the nonlinear term, and the initial data $u_0(\bx) \in L^2(\mathbb{R}^d)$.

\begin{example}\label{example gradient}
\rm
$\mathcal{L} = \vec{b}\cdot\nabla$ is the gradient operator, where $\vec{b}\in \mathbb{R}^d$ is a constant vector. 
\end{example}

\begin{example}\label{example laplace}
\rm
$\mathcal{L} = c \Delta$ is the Laplace operator, where $c>0$ is a constant.
\end{example}

\begin{example}\label{example nonlocal}
	\rm
	$\mathcal{L} = -(-\Delta)^{\alpha / 2}$ is a nonlocal operator with constant $\alpha \in (0,2)$. 
\end{example}

\begin{example}\label{example allencahn}
\rm
For the Allen-Cahn equation, $f(t,\bx,u,\vec{\nabla }u) = u(\bx,t) - (u(\bx,t))^3$.
\end{example}

\begin{example}\label{example hjb}
\rm
For  a special kind of HJB equation, $f(t,\bx,u,\vec{\nabla }u) = \lVert\vec{\nabla }u(\bx,t)\rVert_2^2$,
where $\lVert \cdot \rVert_2$ denotes the standard $L^2$-norm. 
\end{example}

Substituting the change of variable $v(\bx,t) := \me^{-t \mathcal{L}} u(\bx,t)$ into Eq.~\eqref{cauchy problem} yields 
\begin{equation}\label{lawson v}
\frac{\partial}{\partial t} v(\bx,t) = \me^{-t\mathcal{L}} f(t,\bx,\me^{t\mathcal{L}}v,\vec{\nabla }{\left(\me^{t\mathcal{L}}v\right)}),
\end{equation}
which can be readily integrated by any ODE integrator. In view of simplicity and the positivity-preserving property, the first-order time discretization turns out to be more appealing \cite{HochbruckOstermann2010}.  Accordingly, we apply the explicit Euler scheme to Eq.~\eqref{lawson v}:  
\[
v(\bx,t+\tau) \approx v(\bx,t) + \tau \me^{-t\mathcal{L}} f(t,\bx,\me^{t\mathcal{L}}v,\vec{\nabla }{\left(\me^{t\mathcal{L}}v\right)}),
\]
then back to $u(\bx,t)$, we have the so-called Lawson-Euler scheme \cite{lawson1967generalized}:
\begin{equation}\label{Lawson_Euler}
U_{m+1}(\bx) = \me^{\tau \mathcal{L}} \left(U_m(\bx) + \tau  f(t_m,\bx,U_m,\vec{\nabla }U_m)\right),
\end{equation}
where $\tau$ is the time step and $U_m(\bx)$ denotes the numerical solution at $t_m = m\tau$ (we have $U_0(\bx) \equiv u_0(\bx)$ by default). 
Hereafter, we will take the Lawson-Euler scheme \eqref{Lawson_Euler} as an example to show how SPM can be combined with existing time discretization methods.

We use $N$ stochastic particles to evolve Eq.~\eqref{Lawson_Euler}  in the weak sense. The $i$-th particle carries two key quantities: Location $\bx_i \in \mathbb{R}^d$ and weight $w_i \in \mathbb{R}$.  All particles at $t_m$ form a weighted point distribution 
\begin{equation}\label{Xt}
X_{t_m} := \frac{1}{N} \sum_{i=1}^N w_i(t_m) \delta_{\bx_i(t_m)}, 
\end{equation}
where $w_i(t_m)$ and $\bx_i(t_m)$ give the weight and location at time $t_m$ respectively, and $\delta$ is the Dirac distribution. Accordingly, $X_{t_m}$ can approximate $U_m(\bx)$ in the weak sense 
\begin{equation}\label{weak formulation}
\langle \varphi, U_m\rangle \approx \langle \varphi,  X_{t_m} \rangle = \frac{1}{N} \sum_{i=1}^N w_i(t_m) \varphi(\bx_i(t_m)),
\end{equation}
where the standard inner product $\langle f, g\rangle  = \int_{\mathbb{R}^d} f(\bx) g(\bx) \D \bx$ is adopted,   
and the test function $\varphi \in \mathbb{H}^2(\mathbb{R}^d)$ has a compact support in $\mathbb{R}^d$. Various information of $U_m(\bx)$ can be figured out by choosing different test function $\varphi(\bx)$. The weak formulation \eqref{weak formulation} defines a high-dimensional integral that depends on time.

When the particle location $\bx_i$ follows the instrumental distribution $p(\bx)$ at time $t_m$, according to 
\begin{equation}\label{importance sample}
\langle \varphi,  U_m \rangle = \int_{\mathbb{R}^d} \varphi(\bx) \frac{U_m(\bx)}{p(\bx)} p(\bx) \D \bx \approx \frac{1}{N} \sum_{i=1}^{N} \varphi(\bx_i) \frac{U_m(\bx_i)}{p(\bx_i)} = \frac{1}{N} \sum_{i=1}^{N} \varphi(\bx_i) w_i = \langle \varphi,  X_{t_m} \rangle,
\end{equation}
we should let the weight $w_i = {U_m(\bx_i)}/{p(\bx_i)}$. Regardless of what $p(\bx)$ we have chosen, the expectation of $\langle \varphi,  X_{t_m} \rangle$  is always $\langle \varphi, U_m\rangle$, but the variance may be significantly different. The idea of importance sampling \cite{bk:Dimov2008} can be used to reduce the variance. From
\begin{equation}\label{importance sample1}
Var\left(\frac{\varphi U_m}{p}\right) = \int_{\mathbb{R}^d}\left[ \varphi(\bx) \frac{U_m(\bx)}{p(\bx)} \right]^2 p(\bx) \D \bx - (\langle \varphi,  U_m \rangle)^2 = \int_{\mathbb{R}^d}  \frac{(\varphi(\bx)U_m(\bx))^2}{p(\bx)} \D \bx - (\langle \varphi,  U_m \rangle)^2,
\end{equation}
we have the most ideal distribution $p^*(\bx)=\varphi(\bx) U_m(\bx) / \langle \varphi, U_m\rangle$ and thus $Var\left({\varphi U_m}/{p}\right) = 0$
providing that $\varphi(\bx) U_m(\bx)$ is non-negative. However, this is infeasible because: (a) the prerequisite for obtaining $p^*(\bx)$ is that we already know $\langle \varphi, U_m\rangle$ which is our original mission, (b) $\varphi(\bx)$ and $U_m(\bx)$ considered in our work are real-valued functions, and (c) $\varphi(\bx)$ is not predetermined and we use different $\varphi(\bx)$ to extract different information of $U_m(\bx)$. Nevertheless, the idea of importance sampling can still guide us to adjust the distribution of particle location to maintain a small variance. 
More specifically, we are going to use ${|U_m(\bx)|}/{\int_{\mathbb{R}^d}|U_m(\bx)| \D \bx }$  as a feasible alternative to $p^*(\bx)$. At the initial stage, the particle location $\bx_i(0)$ can be sampled from $p(\bx) = {|U_0(\bx)|}/{\int_{\mathbb{R}^d}|U_0(\bx)| \D \bx }$, then assign the weights $w_i(0) = {U_0(\bx_i(0))}/{p(\bx_i(0))}$ accordingly. Afterwards, the Lawson-Euler scheme \eqref{Lawson_Euler} establishes the relationship between $U_m(\bx)$ and $U_{m+1}(\bx)$. It can guide particle motion and weight update. 

In particular, Section \ref{particle motion weight update} gives two evolution strategies of weighted particle system, a preliminary strategy and an improved strategy, denoted as \textbf{Strategy A} and \textbf{Strategy B}, respectively. The derivation of \textbf{Strategy A} is more straightforward, but it poses some restrictions on the nonlinear term $f(t_m,\bx,U_m,\vec{\nabla }U_m)$ while the importance sampling principle is not fully respected. On the contrary, \textbf{Strategy B} can deal with general nonlinear term and respects the importance sampling principle more.

\subsection{ Particle motion and weight update}
\label{particle motion weight update}

To begin with, we introduce \textbf{Strategy A} under the \textbf{Assumption \ref{assum}}.
\begin{assumption}\label{assum}
	\rm
	$U_m(\bx) = 0$ implies $f(t_m,\bx,U_m,\vec{\nabla }U_m) = 0$.
\end{assumption}

From Eq.~\eqref{Lawson_Euler}, we get
\begin{equation}\label{A}
\langle \varphi, U_{m+1} \rangle = \langle \me^{\tau \mathcal{L}^*} \varphi, \left(1 + \tau  \hat{f}(t_m,\bx,U_m,\vec{\nabla }U_m)\right) U_{m} \rangle,
\end{equation}
where $\mathcal{L}^\ast$ denotes the adjoint operator of $\mathcal{L}$ defined on $\mathbb{H}^2(\mathbb{R}^d)$, and 
\begin{equation}
\hat{f}(t_m,\bx,U_m,\vec{\nabla }U_m) :=
\left\{
\begin{aligned}
&\frac{ f(t_m,\bx,U_m,\vec{\nabla }U_m)}{U_m}, &\quad \bx \in \textup{supp}\, U_m, \\
&0, &\quad \bx \notin \textup{supp}\, U_m.
\end{aligned}
\right.
\end{equation}

For \textbf{Example \ref{example allencahn}}, the nonlinear term of Allen-Cahn equation $f(t_m,\bx,U_m,\vec{\nabla }U_m) = U_m(\bx) - \left(U_m(\bx)\right)^3 = 0$ when $U_m(\bx) = 0$ and \textbf{Assumption \ref{assum}} is satisfied. In this case,
\begin{equation}
\hat{f}(t,\bx,U_m(\bx),\vec{\nabla }U_m(\bx)) :=
\left\{
\begin{aligned}
&1-(U_m(\bx))^2, &\quad \bx \in \textup{supp}\, U_m(\bx), \\
&0, &\quad \bx \notin \textup{supp}\, U_m(\bx).
\end{aligned}
\right.
\end{equation}

When $U_m(\bx)$ in Eq.~\eqref{A} is approximated by the weighted point distribution \eqref{Xt}, the dynamics of particle estimator reads that
\begin{equation}\label{dynamics of particle estimator}
\begin{aligned}
 \langle \varphi, U_{m+1} \rangle &= \frac{1}{N} \sum_{i=1}^{N}  w_i(t_m)(1 + \tau  \hat{f}(t_m, \bx_i(t_m),U_m,\vec{\nabla }U_m)) \me^{\tau \mathcal{L}^*} \varphi(\bx_i(t_m))  \\
&\approx \frac{1}{N} \sum_{i=1}^{N} w_i(t_{m+1}) \varphi(\bx_i(t_{m+1})),
\end{aligned}
\end{equation}
where we define
\begin{empheq}[left=\left(\textbf{Strategy A}\right)\empheqlbrace]{align}
	& \bx_i(t_m)\xrightarrow{\me^{\tau \mathcal{L}^*}}\bx_i(t_{m+1}),\label{loca update}\\
	&w_i(t_{m+1}) := w_i(t_m)\left(1 + \tau  \hat{f}(t_m,\bx_i(t_m),U_m,\vec{\nabla }U_m)\right). \label{weight update}
\end{empheq}

In other words, we may directly update the particle weights according to Eq.~\eqref{weight update} where the nonlinear term enters into,
and the dynamics of particle trajectories is updated as Eq.~\eqref{loca update} in the sense of
\begin{equation}\label{x condition expec}
	\mathbb{E}\left[\varphi(\bx_i(t_{m+1})) | \bx_i(t_m)\right] = \me^{\tau \mathcal{L}^*} \varphi(\bx_i(t_m)).
\end{equation}
That is, the undetermined deterministic or stochastic particle dynamics must satisfy Eq.~\eqref{x condition expec} above.
Consider the linear equation,
\begin{equation}\label{adjoint operator linear equation}
\frac{\partial }{\partial t} g (\bx, t)= \mathcal{L}^\ast g(\bx, t), \quad g(\bx, 0) = \varphi(\bx),
\end{equation}
and it has the formal solution $g(\bx,\tau) = \me^{\tau \mathcal{L}^*} \varphi(\bx)$, which reflects the action of $\mathcal{L}^*$ and thus guides the implementation of Eq.~\eqref{loca update}. Some examples are presented below.

For \textbf{Example \ref{example gradient}},  we have $\mathcal{L}^\ast = -\vec{b}\cdot\vec{\nabla}$,
and then the analytical solution of Eq.~\eqref{adjoint operator linear equation} is 
\[
g(\vec{x},\tau) = \varphi(\vec{x}-\vec{b}\tau),
\]
which implies that 
\begin{equation}\label{simu convec}
	\bx_i(t_{m+1}) = \bx_i(t_m) - \vec{b}\tau,
\end{equation}
i.e., a kind of linear advection along the characteristic line.

For \textbf{Example \ref{example laplace}}, we have $\mathcal{L}^\ast = c \Delta$,
and then the analytical solution of Eq.~\eqref{adjoint operator linear equation} has the closed form 
\[
g(\vec{x},\tau) = \int_{\mathbb{R}^d} \frac{1}{(2\pi)^{d/2}\sqrt{|\vec{\Sigma}|}}  \me^{-\frac{1}{2} \vec{y}^T\vec{\vec{\Sigma}}^{-1}\vec{y}} \varphi(\vec{x}-\vec{y})\D \vec{y}, \quad \vec{\Sigma} = 2c \tau \vec{I},
\]
with $\vec{I}$ being the identity matrix of order $d$. It has a transparent probability interpretation $\mathrm{e}^{ c\tau \Delta } \varphi(\vec{x}) = \mathbb{E} \left[\varphi(\vec{x}-\vec{y})\right]$ with $\vec{y}$ obeying the normal distribution $\mathcal{N}(\vec{0}, \vec{\Sigma})$.  To be more specific, it corresponds to the Brownian motion 
\begin{equation}\label{simu diff}
	\vec{x}_i(t_{m+1}) = \vec{x}_i(t_m) - \vec{y}, \;\;\; \text{with} \;\;\; \vec{y} \sim \mathcal{N}(\vec{0}, \vec{\Sigma}).
\end{equation}


The above two examples only involve local operators. For a nonlocal linear operator $\mathcal{L}$, the rules of motions can be designed via the Neumann series expansion. In this  paper, we present a stochastic interpretation of the nonlocal Laplacian operator in \textbf{Example \ref{example nonlocal}}, which has a heat semigroup definition \cite{kwasnicki2017ten}
\begin{equation}\label{fractional heat}
	-(-\Delta)^{\alpha/2} g(\bx, t) = \frac{1}{|\Gamma(-\alpha/2)|} \int_0^{\infty} (\me^{s \Delta} g(\bx, t) - g(\bx, t) ) \frac{\D s}{s^{1 + \alpha/2}},
\end{equation}
where
\begin{equation}
	\me^{s \Delta} g(\bx, t) = \int_{\mathbb{R}^d} \frac{1}{(2\pi)^{d/2}\sqrt{|\vec{\Sigma}|}}  \me^{-\frac{1}{2} \vec{y}^T\vec{\vec{\Sigma}}^{-1}\vec{y}} g(\vec{x}-\vec{y})\D \vec{y}, \quad \vec{\Sigma} =  2s \vec{I}.
\end{equation} 
In accordance with definition \eqref{fractional heat}, we are able to construct a random walk algorithm with power-law tailed jump. \ref{appendix A} presents the details of the random walk algorithm and numerical experiments in solving high-dimensional linear nonlocal equation. 
This kind of processing approach can be generalized to the pseudo-differential operator via Weyl's quantization formula.  For instance, the pseudo-differential operator in the Wigner equation has the probability interpretation of branching random walk \cite{ShaoXiong2020}.

Next, we introduce \textbf{Strategy B} which can handle a general nonlinear term and does not need \textbf{Assumption~\ref{assum}}. From Eq.~\eqref{Lawson_Euler}, we obtain
\begin{equation}
\begin{aligned}
\langle \varphi, U_{m+1} \rangle &= \langle \me^{\tau \mathcal{L}^*} \varphi, U_m + \tau  f(t_m,\bx,U_m,\vec{\nabla }U_m) \rangle \\
&= \left\langle \me^{\tau \mathcal{L}^*} \varphi, \frac{[U_m + \tau  f(t_m,\bx,U_m,\vec{\nabla }U_m)] Z}{|U_m + \tau  f(t_m,\bx,U_m,\vec{\nabla }U_m)|} \cdot \frac{|U_m + \tau  f(t_m,\bx,U_m,\vec{\nabla }U_m)|}{Z} \right\rangle \\
&\approx \frac{1}{N} \sum_{i=1}^{N} \frac{[U_m(\widetilde{\bx_i}(t_m)) + \tau  f(t_m,\widetilde{\bx_i}(t_m),U_m,\vec{\nabla }U_m)] Z}{|U_m(\widetilde{\bx_i}(t_m)) + \tau  f(t_m,\widetilde{\bx_i}(t_m),U_m,\vec{\nabla }U_m)|} \me^{\tau \mathcal{L}^*} \varphi(\widetilde{\bx_i}(t_m)) \\
&\approx \frac{1}{N} \sum_{i=1}^{N} w_i(t_{m+1}) \varphi(\bx_i(t_{m+1})),
\end{aligned}
\end{equation}
where we have applied 
\begin{empheq}[left=\left(\textbf{Strategy B}\right)\empheqlbrace]{align}
	&\bx_i(t_m)\xrightarrow{\text{relocating}} \widetilde{\bx_i}(t_m) \xrightarrow{\me^{\tau \mathcal{L}^*}}\bx_i(t_{m+1}), \label{B 2 motion}\\
	&w_i(t_{m+1}) := \frac{U_m(\widetilde{\bx_i}(t_m)) + \tau  f(t_m,\widetilde{\bx_i}(t_m),U_m,\vec{\nabla }U_m) Z}{|U_m(\widetilde{\bx_i}(t_m)) + \tau  f(t_m,\widetilde{\bx_i}(t_m),U_m,\vec{\nabla }U_m)|} \label{weight update c2},
\end{empheq}
and $\widetilde{\bx_i}(t_m)$ is re-sampled from the density $|U_{m}(\bx) + \tau  f(t_m,\bx,U_m,\vec{\nabla }U_m)|/Z$, 
\begin{equation}\label{x relocating}
\widetilde{\bx_i}(t_m) \sim \frac{1}{Z} |U_{m}(\bx) + \tau  f(t_m,\bx,U_m,\vec{\nabla }U_m)|, \quad Z = \int_{\mathbb{R}^d} |U_{m}(\bx) + \tau  f(t_m,\bx,U_m,\vec{\nabla }U_m)| \D \bx.
\end{equation}
Unlike the location update in Eq.~\eqref{loca update} of \textbf{Strategy A},  
the dynamics of particle trajectories in \textbf{Strategy B} is Eq.~\eqref{B 2 motion} and an extra step of the relocating technique given in Eq.~\eqref{x relocating} is inserted when updating the location from $\bx_i(t_m)$ to $\bx_i(t_{m+1})$. After that, we assign the particle weight according to Eq.~\eqref{weight update c2}. 
\par For \textbf{Strategy A}, the change of particle location is only affected by the linear part as displayed in Eq.~\eqref{loca update}, set against the fact that the evolution of Eq.~\eqref{Lawson_Euler} results from the joint action of linear and nonlinear terms. 
That is, after a period of evolution, the distribution of particle location in \textbf{Strategy A} may not match the shape of solution well (see Figure~\ref{Fig resample reduce variance}), thereby leading to an increase in stochastic variances of Eq.~\eqref{importance sample}. Thanks to the relocating technique in Eq.~\eqref{x relocating}, \textbf{Strategy B} actively adjusts the particle location to obey $|U_{m}(\bx) + \tau  f(t_m,\bx,U_m,\vec{\nabla }U_m)|/Z$, a first-order approximation to ${|U_{m+1}(\bx)|}/{\int_{\mathbb{R}^d}|U_{m+1}(\bx)| \D \bx }$ as shown in the following, 
\begin{equation}
\begin{aligned}
\lvert U_{m+1}(\bx) \rvert &= \lvert \me^{\tau \mathcal{L}} \left(U_m(\bx) + \tau  f(t_m,\bx,U_m,\vec{\nabla }U_m)\right) \rvert \\
&= \left| \left(1+\tau \mathcal{L} + \frac{ \left(\tau \mathcal{L}\right)^2}{2} + \cdots\right) \left(U_m(\bx) + \tau  f(t_m,\bx,U_m,\vec{\nabla }U_m)\right) \right| \\
&= \lvert U_m(\bx) + \tau  f(t_m,\bx,U_m,\vec{\nabla }U_m)\rvert + \mathcal{O}(\tau),
\end{aligned}
\end{equation}
\\
which also agrees with the order of the Lawson-Euler scheme~\eqref{Lawson_Euler}. When the evolution is performed in a short time interval, the difference between the two strategies may be not significant. But in long-time simulations, the relocating technique in Eq.~\eqref{x relocating}, making the particle location close to the shape of $U_m(\bx)$, is crucial, with which the particle weight is more uniform and the variance is smaller when calculating Eq.~\eqref{importance sample}. This can be readily verified through benchmark tests in Section~\ref{sec:acc}. 

\subsection{Piecewise constant reconstruction from weighted point distribution}
\label{section approximation nonlinear}

According to the dynamics of particle estimator in Eqs.~\eqref{weight update} and \eqref{weight update c2}, we need to know the values of $U_m(\bx)$ and $f(t,\bx,U_m(\bx),\vec{\nabla }U_m(\bx))$ for both strategies.
It requires to approximate $U_m(\bx)$, $\vec{\nabla }U_m(\bx)$ using the weighted point distribution \eqref{Xt}. 
Generally speaking, when $U_m(\bx)$ is sufficiently smooth, those quantities like  $\vec{\nabla }U_m(\bx)$ can be approximated by corresponding difference quotients, so we just need to focus on how to compute $U_m(\bx)$.  There exist several numerical techniques for constructing function values from a point cloud, such as the piecewise constant approximation \cite{raviart1983particle}, the Fourier series approximation \cite{Yan2015} and the sequential function approximation \cite{wu2017randomized, wu2018sequential}. 
In this work, we adopt the first-order piecewise constant approximation to obtain $U_m(\bx)$.

Let $\Omega = \prod\limits_{j=1}^d \left[l_j,r_j\right] \subset \mathbb{R}^d$ be the computational domain containing all particles.  Suppose the particle location $\bx_i(t_m) = \left(x_i^1,x_i^2,\dots,x_i^d\right),\, i = 1,2,\dots,N$, and let $l_j = \min\limits_{1\leq i\leq N} x_i^j$ and $r_j = \max\limits_{1\leq i\leq N} x_i^j$. Then, the computational domain can be decomposed into
\begin{equation}\label{eq:decom}
\Omega={\prod\limits_{j=1}^d \left[l_j,r_j\right]} = \bigcup_{k=1}^K Q_k,
\end{equation}
so that $U_m(\bx)$ can be approximated by constant in each sub-domain $Q_k$ (usually a hypercube with the size length of $h_k$). The discussion on how to achieve such decomposition is left to Section \ref{implement section}. Suppose location $\vec{x}$ at which the function value needs to be evaluated falls into the sub-domain $Q_k$, then in a similar way to Eq.~\eqref{weak formulation}, we have
\begin{equation}\label{eq weight sum approx u}
U_m(\bx) \approx  \frac{1}{\mu(Q_k)}\int_{Q_k} U_m(\by) \D \vec{y} \approx \frac{1}{N}\sum_{i=1}^N \frac{w_i(t_m) \mone_{Q_{k}}(\bx_i(t_m))}{\mu(Q_k)},
\end{equation}
where $\mu(Q_k)$ gives the Lebesgue measure of $Q_k$. In fact, we are able to obtain the second approximation in Eq.~\eqref{eq weight sum approx u} 
by taking $\varphi(\bx)$ as the indicator function $\mone_{Q_{k}}(\bx)$ in Eq. \eqref{weak formulation}. Let $\bx_k^c$ be the center point of $Q_k$. If $\bx = \bx_k^c$, then $U_m(\bx) = \frac{1}{\mu(Q_k)}\int_{Q_k} U_m(\by) \D \vec{y} + \mathcal{O}(h_k^2)$, otherwise $U_m(\bx) = \frac{1}{\mu(Q_k)}\int_{Q_k} U_m(\by) \D \vec{y} + \mathcal{O}(h_k)$.

The nonlinear term in \textbf{Example \ref{example hjb}} contains a gradient operator and we adopt the central difference quotient to approximate the derivative in the $j$-th coordinate,  
\begin{equation}\label{center quotient}
	\begin{aligned}
		\frac{\partial U_m(\bx)}{\partial x_j} &\approx \frac{\partial U_m(\bx_k^c)}{\partial x_j} \approx  \frac{U_m(\bx_k^c+h_k\be_j)-U_m(\bx_k^c-h_k\be_j)}{2 h_k}\\
		&\approx \frac{\frac{1}{\mu(Q_{k_r})}\int_{Q_{k_r}} U_m(\by) \D \vec{y} -\frac{1}{\mu(Q_{k_l})}\int_{Q_{k_l}} U_m(\by) \D \vec{y}}{2 h_k},
	\end{aligned}
\end{equation}
and then using Eq.~\eqref{eq weight sum approx u} yields the approximation to $\frac{\partial U_m(\bx)}{\partial x_j}$, 
where $\bx_k^c+h_k\be_j$ and $\bx_k^c-h_k\be_j$ are supposed to fall into the sub-domain $Q_{k_r}$ and $Q_{k_l}$, respectively, and the unit vector $\be_j$ denotes the direction of the $j$-th coordinate. In actual numerical simulations, we choose $\bx_k^c+h_k\be_j$ and $\bx_k^c-h_k\be_j$ to respectively be the center of $Q_{k_r}$ and $Q_{k_l}$ for ensuring that the approximation error in Eq.~\eqref{center quotient} is $\mathcal{O}(h_k)$.


Up to this point, we have delineated the mathematical description of SPM and are also able to classify its errors into three levels,
\begin{itemize}
\item[$\bullet$] Spatial discretization: The piecewise constant reconstruction~\eqref{eq weight sum approx u} introduces a spatial approximation error of $\mathcal{O}(h_k)$;
\item[$\bullet$] Time discretization: The Lawson-Euler scheme~\eqref{Lawson_Euler} introduces a time discretization error of $\mathcal{O}(\tau)$;
\item[$\bullet$] Monte Carlo sampling: The stochastic particle approximation~\eqref{weak formulation} introduces a sampling error of $\mathcal{O}(N^{-1/2})$.
\end{itemize}
Generally speaking, the sampling error dominates, especially for the solutions with complex landscapes. For the initial sampling, there is a significant body of research on efficiently sampling from complex distributions \cite{Roberts2001metropolis}.
{We would like to remark on other kinds of time discretization rather than the Lawson-Euler scheme~\eqref{Lawson_Euler} including, for example,  the exponential time differencing scheme \cite{beylkin1998new}
and the Strang splitting scheme \cite{hansen2012second}. 
The Lawson-Euler scheme may be more stable than the exponential time differencing scheme \cite{crouseilles2020exponential}, though they are both of the first order accuracy. Besides, considering the time discretization error does not dominate the errors of SPM in general, 
even the first-order accurate Lawson-Euler scheme may allow a large time step without recourse to the second-order accurate Strang splitting.}

\section{Implementation details and cost analysis}
\label{implement section}

We have already had two strategies for particle motion and weight update, but there still remain some implementation details, such as how to achieve the decomposition \eqref{eq:decom} and how to sample $\widetilde{\bx_i}(t_m)$ in \eqref{x relocating}.

\subsection{Piecewise constant reconstruction and virtual uniform grid (VUG)}
\label{virtual uniform grid subsection}

\begin{figure}[htbp]
	\centering
	\includegraphics[width=0.5\linewidth]{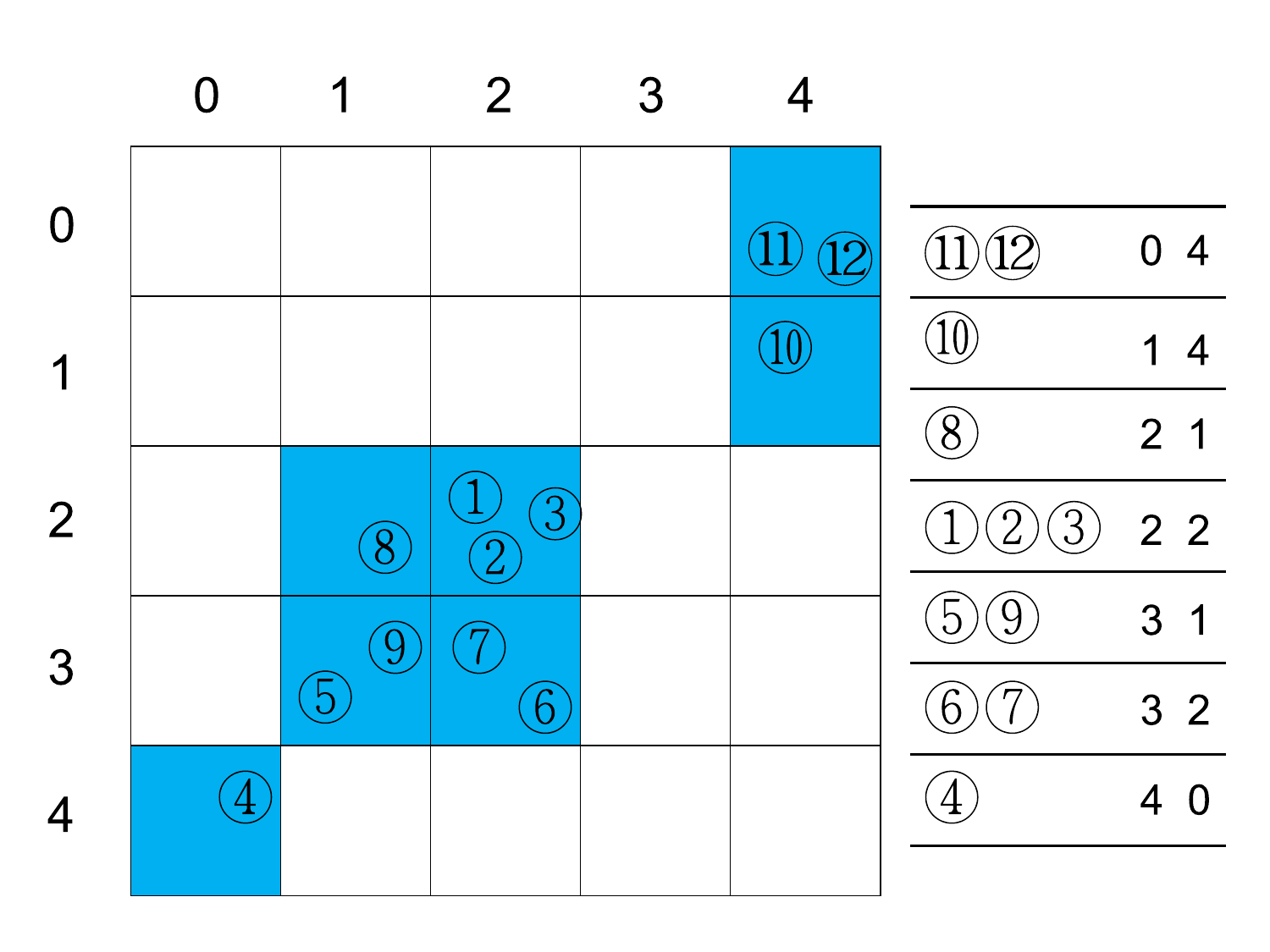}
	\caption{VUG: Only blue grids, containing the particles, are stored, and the corresponding grid coordinates for all particles are presented on the right.}
	\label{vug schematic diagram}
\end{figure}

 An optimal way for piecewise constant reconstruction can be implemented by solving  
\begin{equation} \label{optimal problem}
	\inf\limits_{Q_k} \sup\limits_{\varphi \in D(\varphi)} \left| \int_{\mathbb{R}^d} \varphi(\bx) \sum\limits_{Q_k} \left(\sum_{i=1}^{N} \frac{w_i \mone_{Q_{k}}(\bx_i(t_m))}{ N\mu\left({Q_k}\right)}\right) \mone_{Q_{k}}(\bx)  \D \bx - \frac{1}{N} \sum_{i=1}^N w_i \varphi(\bx_i(t_m))  \right|,
\end{equation}
where $D(\varphi) = \{\varphi \in  \mathbb{H}^2(\mathbb{R}^d): \varphi \text{ has bounded variation}\}$. 
Along this line, an adaptive sequential partition strategy based on the star discrepancy, which measures the degree of uniformity, was proposed in \cite{LiYangWong2016}, where the decomposition is sequentially refined until the star discrepancy of particle locations inside each $Q_k$ is less than a prescribed threshold. 
However, the complexity of calculating the star discrepancy is NP-hard \cite{GnewuchWahlstromWinzen2012} and  
one has to deal with a significant computational cost when $N$ is large. 
Instead, considering the highly adaptive nature of stochastic particles, this work chooses a ``uniform partitioning" as an approximation to the problem \eqref{optimal problem} without storing the full {tensor} girds, and the resulting method is termed ``virtual uniform grid" (VUG). {More specifically, during the $m$-th step, the particle location is close to ${|U_m(\bx)|}/{\int_{\mathbb{R}^d}|U_m(\bx)| \D \bx }$ and VUG does not store the grids inside areas where ${|U_m(\bx)|}/{\int_{\mathbb{R}^d}|U_m(\bx)| \D \bx }$ vanishes, thereby avoiding unnecessary storage as much as possible.}
Given a uniform decomposition as shown in Eq.~\eqref{eq:decom}, each $Q_k$ is a hypercube with the same side length, denoted by $h$,
and thus VUG has a storage cost of around $\alpha \prod\limits_{j=1}^{d}\left(\frac{r_j-l_j}{h}\right)$, 
where $\alpha\in (0,1)$ is the ratio of the stored grids to the full ones.  Figure~\ref{vug schematic diagram}
gives an illustration where only 7 grids containing particles in all 25 grids should be stored, namely $\alpha = \frac{7}{25}$. In general, such ratio $\alpha$ depends on the shape of unknown functions, and generally decreases as the dimension increases, which can be readily verified,  for instance, through the numerical experiments in \textbf{Example~\ref{exm:vug}}.


\begin{example}\label{exm:vug}
	\rm 	
	Consider {a} sign-changing function 
	\begin{equation}
		p^\text{ref}(\bx) = \prod_{j=1}^d f_{15, 5}(x_j) - \prod_{j=1}^d f_{10, 10}(x_j) + \prod_{j=1}^d f_{5, 15}(x_j)
	\end{equation}
	with $f_{\alpha, \beta}(x) $ being the beta distribution density
	\begin{equation}
		f_{\alpha, \beta}(x) = \frac{\Gamma(\alpha + \beta)}{\Gamma(\alpha) \Gamma(\beta)} x^{\alpha - 1}(1 - x)^{\beta -1}, \quad 0 \le x \le 1.
	\end{equation}
	One can sample from $p^\text{ref}(\bx)$ and do piecewise constant reconstruction through VUG to get a numerical solution $p^\text{num}(\bx)$. Table~\ref{VUGtable} shows the numerical results of VUG in different dimensions $d$ for different sample sizes $N$ where the relative $L^2$ error 
	\begin{equation}\label{L2 static}
		\mathcal{E}_2[p] = \frac{\lVert p^\text{num}(\bx)-p^\text{ref}(\bx)\rVert_2}{\lVert p^\text{ref}(\bx)\rVert_2}
	\end{equation} 
	is adopted to investigate the accuracy. It can be easily observed there that  $\mathcal{E}_2[p]$ decreases as the sample size $N$ increases for fixed $d$, and the ratio $\alpha$ always decreases as the dimension $d$ increases.

\end{example}

\begin{table}[htbp]
	\centering
	\caption{Example~\ref{exm:vug}: The relative $L^2$ error and the ratio $\alpha$ when using VUG to do piecewise constant reconstruction in 4-D to 7-D with a fixed side length $h=0.0625$. \#~VUG in the fourth column represents the number of virtual uniform grids, {$\prod\limits_{j=1}^{d}\left(\frac{r_j-l_j}{h}\right)$} in the fifth column represents the number of full grids. Overall, the ratio $\alpha$ decreases as the dimension $d$ increases.}
	\label{VUGtable}
	\begin{tabular}{cccccc}
		\toprule
		$d$ & $N$       &  {$\mathcal{E}_2[p]$}  & \#~VUG   &  {$\prod\limits_{j=1}^{d}\left(\frac{r_j-l_j}{h}\right)$} & $\alpha$   \\ \midrule
		& $1\times 10^6$ &0.0693 & $1.86\times 10^4$&  & 0.283\\
		4 & $2\times 10^6$ &0.0573 & $2.09\times 10^4$& $6.55\times 10^4$ & 0.318\\
		& $4\times 10^6$ &0.0495 & $2.31\times 10^4$&  & 0.353\\ \midrule
		& $1\times 10^7$ &0.0657 & $1.75\times 10^5$&  & 0.166\\
		5 & $2\times 10^7$ &0.0584 & $2.01\times 10^5$& $1.05\times 10^6$ & 0.191\\
		& $4\times 10^7$ &0.0543 & $2.29\times 10^5$&  & 0.218\\ \midrule
		& $4\times 10^7$ & 0.0777 &  $1.26\times 10^6$   &  & 0.075  \\ 
		6   & $8\times 10^7$ & 0.0691 &$1.52\times 10^6$  &$1.68\times 10^7$ & 0.091  \\ 
		& $1.6\times 10^8$  & 0.0641 &$1.80\times 10^6$  &  & 0.107    \\ \midrule
		& $1\times 10^8$ & 0.1001 & $7.00\times 10^6$ &   & 0.026  \\ 
		7   & $2\times 10^8$ & 0.0861 &$8.94\times 10^6$  &$2.68\times 10^8$ & 0.033 \\ 
		& $4\times 10^8$  & 0.0776 & $1.12\times 10^7$ &    & 0.041  \\
		\bottomrule
	\end{tabular}
	
\end{table}

\begin{figure}[htbp]
	\centering
	\includegraphics[width=0.5\textwidth]{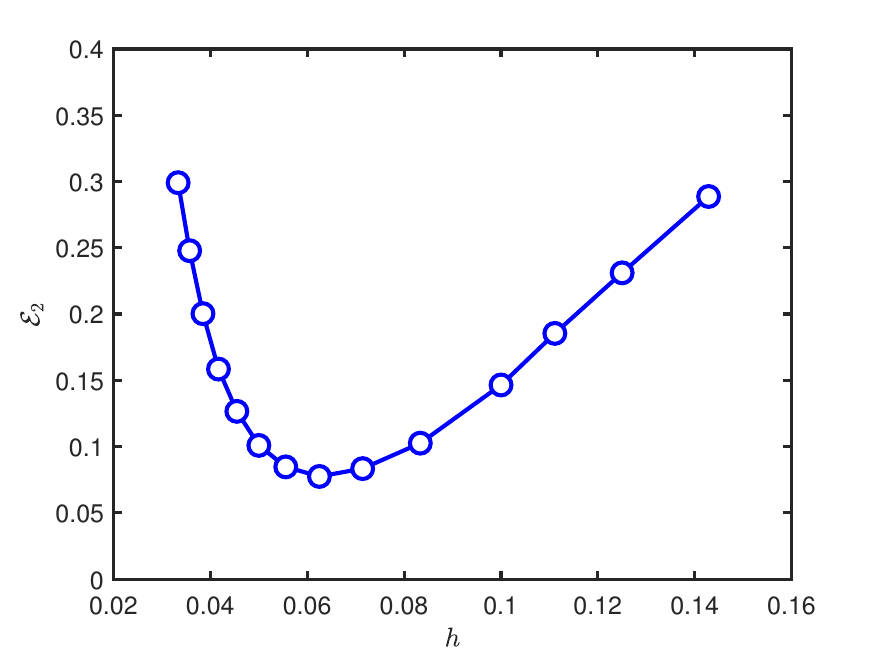}
	\caption{Example~\ref{exm:vug}: The relative $L^2$ error $\mathcal{E}_2[p]$ against the  side length $h$ for $d=6$ and $N = 4\times 10^7$.}
	\label{table density err h}
\end{figure}

The uniform side length $h$ of hypercubes in the decomposition \eqref{eq:decom} introduces 
an asymptotic error of $\mathcal{O}(h)$ for the piecewise constant reconstruction through VUG \cite{XiongShao2019},
and thus a large $h$ is not preferred. For \textbf{Example~\ref{exm:vug}}, such asymptotic error becomes more evident for 
$h > 0.1$, see  Figure~\ref{table density err h}. However, too small $h$ should NOT be suggested due to the  notorious overfitting problem \cite{Yan2015,bk:Silverman2018}, which results in a significant statistical error for the number of grids exceeds the sample size and few particles fall within the same grid. This is also verified in Figure~\ref{table density err h}:  When $h<0.04$, the effective sample size in a grid becomes extremely small and leads to an augmentation of statistical error. A practical way to choose an appropriate $h$ may be as follows. We may take a relatively large $h$ first, and then gradually reduce it until the relative $L^2$  error $\mathcal{E}_2[p]$ does not diminish with the decrease of $h$.

In practice, VUG dynamically allocates the memory according to the particle locations at each time step, whereas the memory is usually allocated in advance for the full tensor grids. Within our C++ implementation, the data structure  ``map" is adopted to map the grid coordinates in struct type to piecewise constants in double type. Specifically, we need to maintain two ``map" data structures: $\mathcal{M}_1$ records the piecewise constant approximation to the unknown function and $\mathcal{M}_2$ the piecewise constant approximation to the nonlinear term. Algorithm \ref{alg SPM} presents the pseudo-code of  SPM with VUG and Algorithm \ref{alg m1} details the updating $\mathcal{M}_1$  for an incoming new particle.

\begin{algorithm}[htbp]
	\caption{ Stochastic particle method with virtual uniform grid. } 
	\label{alg SPM}
	\hspace*{0.02in} {\bf Input:}
	The computational domain ${\prod\limits_{j=1}^d \left[l_j,r_j\right]}$, the side length of hypercube $h$, 
	the final time $T$, the time step $\tau$, the initial data $u_0(\bx)$, the linear operator $\mathcal{L}$ and the nonlinear term {$f(t,\bx,u,\vec{\nabla }u)$}.\\
	\hspace*{0.02in} {\bf Output:} The numerical solution to nonlinear PDE \eqref{cauchy problem} at the final time.
	\begin{algorithmic}[1]
		\State Allocate two empty maps $\mathcal{M}_1$ and $\mathcal{M}_2$;
		\For{$m=0:\frac{T-\tau}{\tau}$} 
		\State Clear $\mathcal{M}_1$;
		\If{\textbf{Strategy A}}
		\For{$i=1:N$}
		\If{$m = 0$}
		\State Sample $\bx_i$ from $p(\bx)=\frac{|U_0(\bx)|}{\int_{\mathbb{R}^d} |U_0(\bx)| \D \bx}$;
		 \State $w_i\leftarrow\frac{U_0(\bx_i)}{{p(\bx_i)}}$;
		\EndIf
		\State Update $w_i$ according to \eqref{weight update};
		\State Move particle according to $\mathcal{L}^*$ (the adjoint operator of $\mathcal{L}$);
		\State Update $\mathcal{M}_1$ according to $\bx_i$ and $w_i$;
		\EndFor
		\State Clear $\mathcal{M}_2$;
		\State  Compute $\hat{f}(t_{m+1},\bx,U_{m+1},\vec{\nabla }U_{m+1 })$ and construct $\mathcal{M}_2$ based on $\mathcal{M}_1$; 
		\ElsIf{\textbf{Strategy B}}
		\For{${i=1:N}$}
		\State Sample $\bx_i$ from ${p(\bx)}=\frac{|U_{m}(\bx) + \tau  f(t_m,\bx,U_m,\vec{\nabla }U_m)|}{\int_{\mathbb{R}^d} |U_{m}(\bx) + \tau  f(t_m,\bx,U_m,\vec{\nabla }U_m)| \D \bx}$;
		\State  $w_i\leftarrow\frac{U_{m}(\bx_i) + \tau  f(t_{m},\bx_i,U_{m},\vec{\nabla }U_{m})}{{p(\bx_i)}}$;
		\State Move particle according to $\mathcal{L}^*$;
		\State Update $\mathcal{M}_1$ according to $\bx_i$ and $w_i$;
		\EndFor
		\State Clear $\mathcal{M}_2$;
		\State Compute $U_{m+1}(\bx) + \tau  f(t_{m+1},\bx,U_{m+1},\vec{\nabla }U_{m+1})$ and construct $\mathcal{M}_2$ based on $\mathcal{M}_1$;
		\State Compute $\int_{\mathbb{R}^d} |U_{m+1}(\bx) + \tau  f(t_{m+1},\bx,U_{m+1},\vec{\nabla }U_{m+1})| \D \bx$;
		\EndIf
		\EndFor
	\end{algorithmic}
\end{algorithm}

\begin{algorithm}[htbp]
	\caption{Update {piecewise constant approximation} $\mathcal{M}_1$ according to a new particle.} 
	\hspace*{0.02in} {\bf Input:}
	\label{alg m1}
	An not updated $\mathcal{M}_1$, a new particle with location $\bx = \left(x_1,x_2,\dots,x_d\right)$ and weight $w$,
	the computational domain ${\prod\limits_{j=1}^d \left[l_j,r_j\right]}$ and the side length of hypercube $h$. \\
	\hspace*{0.02in} {\bf Output:} The updated $\mathcal{M}_1$.
	\begin{algorithmic}[1]
		\State Compute the grid coordinate {$\mathbf{coord} = \left(\text{floor}\left(\frac{x_1 - l_1}{h}\right), \text{floor}\left(\frac{x_2 - l_2}{h}\right),\dots,\text{floor}\left(\frac{x_d - l_d}{h}\right)\right)$;}
		\If{$\mathbf{coord}$ already exists in $\mathcal{M}_1$}
		\State $\mathcal{M}_1(\mathbf{coord}) \leftarrow \mathcal{M}_1(\mathbf{coord}) + \frac{w}{N h^d}$;
		\Else
		\State Insert $\mathbf{coord}$ into $\mathcal{M}_1$;
		\State $\mathcal{M}_1(\mathbf{coord}) \leftarrow  \frac{w}{N h^d}$;
		\EndIf
	\end{algorithmic}
\end{algorithm}

Suppose we have obtained the piecewise constant reconstruction 
of $U_m(\bx)$ via VUG and thus $U_{m}(\bx) + \tau  f(t_m,\bx,U_m,\vec{\nabla }U_m) = \sum_{k=1}^K c_k \cdot \mathbbm{1}_{Q_k}(x)$ is also a piecewise constant function,  then we have 
\begin{equation}\label{compute Z}
	Z = \int_{\mathbb{R}^d} |U_{m}(\bx) + \tau  f(t_m,\bx,U_m,\vec{\nabla }U_m)| \D \bx  = \sum_{k=1}^K |c_k| \cdot h^d.
\end{equation}
To compute $Z$, it is only needed to sum the elements in $\mathcal{M}_2$ and multiply the sum by $h^d$. In order to implement \textbf{Strategy B}, we need to further sample the particle location $\bx_i$ from the piecewise constant density distribution $|U_{m}(\bx) + \tau  f(t_m,\bx,U_m,\vec{\nabla }U_m)|/Z$. {The probability that $\bx_i$ falls in $Q_k$ is ${|c_k|h^d}/{Z}$, and it is uniformly sampled within $Q_k$.}

\subsection{Domain decomposition and MPI parallelization} 
\label{sec mpi}

The distributed computing technology via MPI as well as some domain decomposition strategies are adopted to parallelize SPM and  the pseudo-code is presented by Algorithm~\ref{domain decomposition alg}.  The computational domain $\Omega$ in the uniform decomposition Eq.~\eqref{eq:decom} is decomposed into $nproc$  blocks, with $nproc$ being the number of MPI processes and usually taken as a power of two for convenience. Each block is a union of some hypercubes $Q_k$ with the same side length $h$. At the beginning of domain decomposition, we let a domain set $\mathcal{S} = \{\Omega\}$, and each domain in $\mathcal{S}$ is sequentially divided into two domains (still use $\mathcal{S}$ to collect the resulting domains) until $|\mathcal{S}|=nproc$. Suppose $\mathcal{S} \ni \Omega_{s} = \prod \limits_{j=1}^d \left[a_j^{(s)},b_j^{(s)}\right] $
and $c_{i_0,j_0}^{(s)}=a_{j_0}^{(s)}+i_0\times h$ are some partition nodes to be selected
where $i_0 = 1,2,\dots,{(b_{j_0}^{(s)}-a_{j_0}^{(s)})}/{h}-1$, $j_0 = 1,2,\dots,d$. To strike the load balance, it is {expected} that the number of particles in each block is roughly the same. To this end, 
the optimal {split} points, denoted by $\widetilde{c}_{i_0,j_0}^{(s)}$,   are required to minimize the difference
\begin{equation}\label{domain decom gjk}
\widetilde{c}_{i_0,j_0}^{(s)} = 
\argmin_{c_{i_0,j_0}^{(s)} }  | L_{i_0,j_0}^{(s)} - R_{i_0,j_0}^{(s)}|,
\end{equation}
where $\Omega_{s}$ is divided into $\Omega_{s}^{(1)}$ and $\Omega_{s}^{(2)}$: 
\begin{align}
	\Omega_{s}^{(1)} &= \prod \limits_{j=1}^{j_0-1} \left[a_j^{(s)},b_j^{(s)}\right] \times [a_{j_0}^{(s)},\widetilde{c}_{i_0,j_0}^{(s)}]\times  \prod \limits_{j=j_0+1}^{d}  \left[a_j^{(s)},b_j^{(s)}\right],  \label{omega1} \\ 
	\Omega_{s}^{(2)} &= \Omega_s \backslash \Omega_s^{(1)}, \label{omega2}
\end{align}
and $L_{i_0,j_0}^{(s)}$ and $R_{i_0,j_0}^{(s)}$ count the number of particles that fall in $\Omega_s^{(1)}$ and $\Omega_s^{(2)}$, respectively.  In order to compute $L_{i_0,j_0}^{(s)}$ and $R_{i_0,j_0}^{(s)}$ {efficiently}, we only consider the first two coordinate components and restrict $j_0\leq 2$. To be more specific, we create a 2-D array $A$ of size {$\frac{r_1-l_1}{h} \times \frac{r_2-l_2}{h}$}, and $A[i][j]$ records the number of particles whose the first two location components fall in the {$\left[l_1+ih,l_1+(i+1)h\right) \times \left[l_2+jh,l_2+(j+1)h\right)$}. In this way, $L_{i_0,j_0}^{(s)}$ and $R_{i_0,j_0}^{(s)}$ can be readily obtained according to the 2-D array $A$ and thus the domain decomposition can be done conveniently.

\begin{algorithm}[htbp]
\caption{Domain decomposition.} 
\label{domain decomposition alg}
\hspace*{0.02in} {\bf Input:} 
 The 2-D array $A$, the computational domain $\Omega$, the side length of hypercube $h$ and the process number $nproc$. \\
\hspace*{0.02in} {\bf Output:} 
The domain set $\mathcal{S}=\{\Omega_1,\Omega_2,\dots,\Omega_{nproc}\}$.
\begin{algorithmic}[1]
\State Let $\mathcal{S}$ {an empty} domain set;
\State Let $\mathcal{B}$ {an empty} queue (first in, first out) of domain, push $\Omega$ into $\mathcal{B}$;
\State $n = 1$;

\While{$n < nproc$} 
\State {Pick up $\Omega_s$ as} {the top} element of $\mathcal{B}$;
\State Pop top element of $\mathcal{B}$;
\State Compute $L_{i_0,j_0}^{(s)}$ and $R_{i_0,j_0}^{(s)}$ according to $\Omega_s$ and $A$;
\State Choose $\widetilde{c}_{i_0,j_0}^{(s)}$ to attain the minimum of Eq.~\eqref{domain decom gjk};
\State Divide domain $\Omega_s$ into $\Omega_s^{(1)}$ and $\Omega_s^{(2)}$ according to Eqs.~\eqref{omega1} and \eqref{omega2};
\State Push $\Omega_s^{(1)}$ and $\Omega_s^{(2)}$ into $\mathcal{B}$;
\State $n \leftarrow n+1$;
\EndWhile
\State Put all elements of $\mathcal{B}$ into domain set $\mathcal{S}$; 
\end{algorithmic}
\end{algorithm}

\subsection{Cost analysis}
\begin{itemize}

\item[$\bullet$] Storage cost: 
{Storing} virtual uniform grids costs {$\mathcal{O}(\alpha \prod\limits_{j=1}^{d}\left(\frac{r_j-l_j}{h}\right))$} for both \textbf{Strategy A} and \textbf{Strategy B}. In practice, the number of grids is set to not exceed the number of particles, therefore the memory cost is $\mathcal{O}(N)$.  In addition, \textbf{Strategy A} needs an extra memory to store $N$ particles because $\bx_i(t_m)$ and $w_i(t_m)$ are required at $t_{m+1}$ according to Eqs.~\eqref{loca update} and \eqref{weight update}. 

\item[$\bullet$] Time complexity: 
The most time-consuming parts are query and insertion of map required in the piecewise constant reconstruction as shown in Algorithm \ref{alg m1}, and its complexity is bounded by $\mathcal{O}(N\log N)$ since {the number of grids never exceeds the number of particles.}
For the linear operators in \textbf{Example \ref{example gradient}} and \textbf{Example \ref{example laplace}},
{the complexity of particle motion is $\mathcal{O}(N)$.} That is, the complexity of \textbf{Strategy A} is $\mathcal{O}(N\log N)$.
Moreover, \textbf{Strategy B} performs two more $\mathcal{O}(N)$ operations than \textbf{Strategy A}: Relocating in Eq.~\eqref{x relocating} and computing $Z$ in Eq.~\eqref{compute Z}. 

\end{itemize}

\section{Accuracy test}
\label{sec:acc}

\begin{figure}[H]
	\centering
	\begin{subfigure}[b]{0.49\textwidth}
		\centering
		\includegraphics[width=\textwidth]{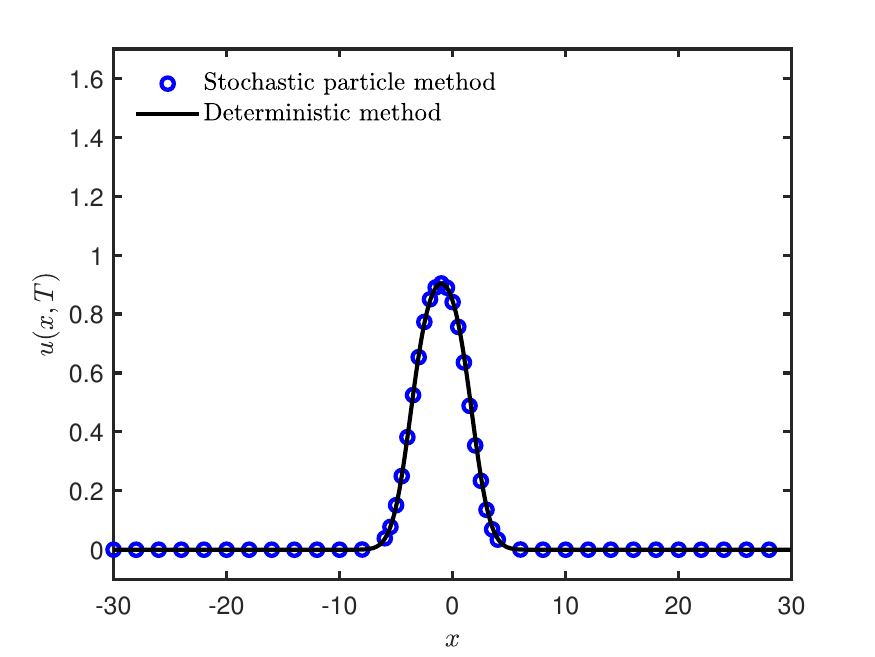}
		\caption{$T=1$, \textbf{Strategy A}.}
	\end{subfigure}
	\hfill
	\begin{subfigure}[b]{0.49\textwidth}
		\centering
		\includegraphics[width=\textwidth]{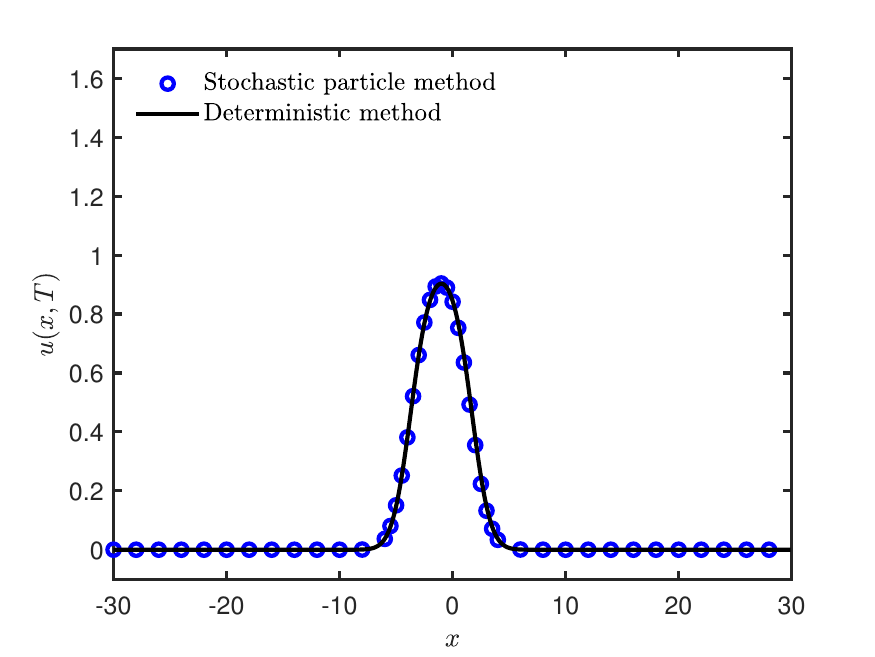}
		\caption{$T=1$, \textbf{Strategy B}.}
	\end{subfigure}
	\begin{subfigure}[b]{0.49\textwidth}
		\centering
		\includegraphics[width=\textwidth]{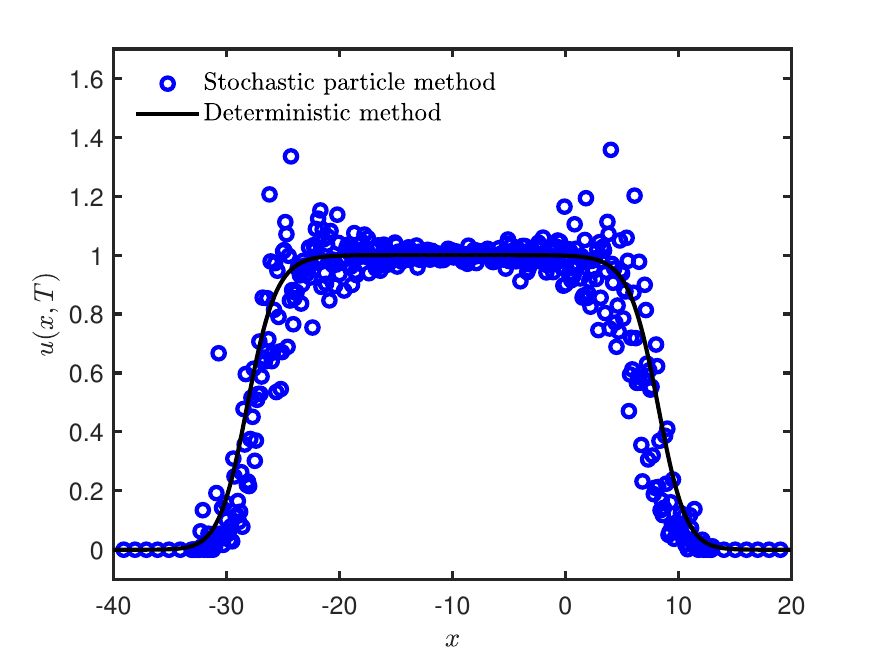}
		\caption{$T = 10$, \textbf{Strategy A}.}
	\end{subfigure}
	\hfill
	\begin{subfigure}[b]{0.49\textwidth}
		\centering
		\includegraphics[width=\textwidth]{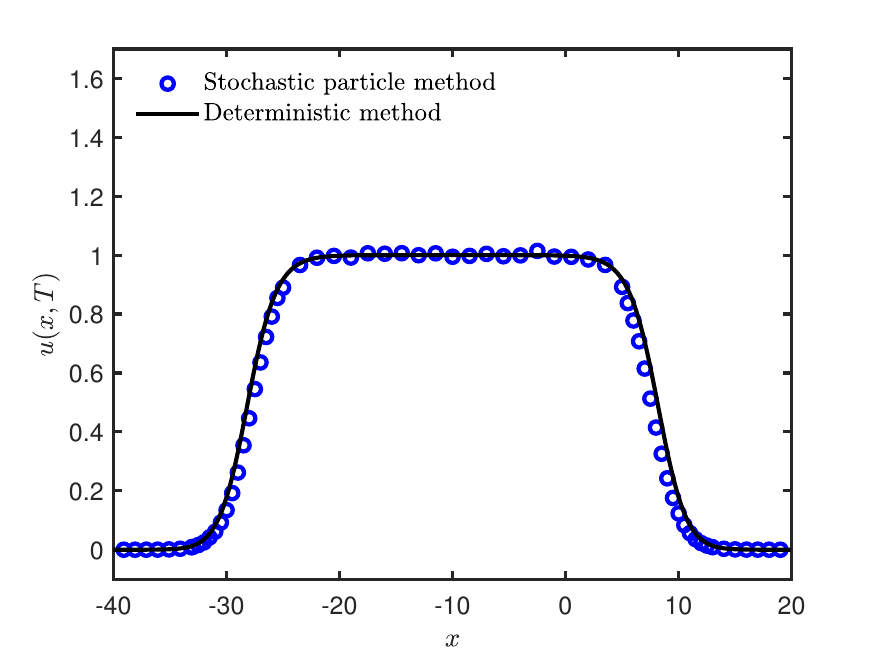}
		\caption{$T=10$, \textbf{Strategy B}.}
	\end{subfigure}
	\caption{\textbf{Strategy A} (\textbf{preliminary strategy}) v.s. \textbf{Strategy B} (\textbf{improved strategy}): The relocating technique (see Eq.~\eqref{x relocating}) adopted by  \textbf{Strategy B} may reduce the variance in long-time simulations. Here we set $N=1.6\times 10^7$, $\tau=0.1$ and $h=0.1$.}
	\label{Fig resample reduce variance}
\end{figure}

In order to benchmark SPM, we consider the following 1-D nonlinear model with a Cauchy data
\begin{equation}\label{eq:1DAC}
u_t - { \nabla u} - \Delta u - u + u^3 = 0, \quad u(x,0) = u_0(x) = \exp(-x^2)(1+x^4),
\end{equation}
and employ the relative $L^2$ error 
\begin{equation}\label{eq:L2}
\mathcal{E}_2[u](t) = \frac{\lVert u^\text{num}(x, t)-u^\text{ref}(x, t)\rVert_2}{\lVert u^\text{ref}(x, t)\rVert_2},
\end{equation}
to measure the accuracy. The reference solutions $u^\text{ref}(x, t)$ to Eq.~\eqref{eq:1DAC} are produced by a deterministic solver which adopts the second-order operator splitting with two exact flows \cite{hansen2012second}. {The nonlinear terms {$f(t,x,u,\nabla u) = u-u^3$} in Eq.~\eqref{eq:1DAC} satisfy \textbf{Assumption \ref{assum}}, thus both \textbf{Strategy A} and \textbf{B} can be tested.} As expected,  SPM {shows a half-order convergence} with respect to the sample size $N$,
and a first-order accuracy in $\tau$ and $h$ when \textbf{Strategy A} and \textbf{Strategy B} are implemented with the final time $T=1$ and $T=10$, respectively (see Table~\ref{tableABaccuracy}).
However, the numerical solution of \textbf{Strategy A} at $T=10$ deteriorates 
due to the accumulated stochastic variance while that of \textbf{Strategy B} does not (see Figure~\ref{Fig resample reduce variance}). This exactly reflects the effectiveness of the relocating technique in Eq.~\eqref{x relocating} adopted by \textbf{Strategy B} for long-time simulations. Actually, as shown in Figure~\ref{Fig resample reduce variance}, the difference of \textbf{Strategy A} and \textbf{Strategy B} is not apparent until $T=1$ because the particle location matches the shape of solution (close to the initial shape of {$u_0(x)$}) well for both strategies. 
At $T=10$, {the shape of solutions is evidently changed and} the distribution of particle location in \textbf{Strategy A} may not match it without relocating. By contrast, the relocating technique given in Eq.~\eqref{x relocating} helps \textbf{Strategy B} to actively {adjust} the particle location and thus match the shape of solution within the allowed error tolerance, thereby reducing the stochastic variance of Eq.~\eqref{importance sample}. 
Therefore, the implementation of SPM with \textbf{Strategy B} will be used to update the particle system in 
Sections \ref{sec:ac} and \ref{sec:hjb} for solving moderately high-dimensional problems. 

\bigskip

\begin{table}[H]
	\centering
	\caption{Convergence of SPM: First-order accuracy in both time and space plus half-order accuracy in the sample size. \textbf{Strategy A} and \textbf{Strategy B} are implemented with the final time $T=1$ and $T=10$, respectively. }
	\begin{tabular}{ccccccc}
		\hline
		\multirow{2}{*}{$N$} & \multirow{2}{*}{$h$} & \multirow{2}{*}{$\tau$} &\multicolumn{2}{c}{\textbf{Strategy A}, $T = 1$} &\multicolumn{2}{c}{\textbf{Strategy B}, $T = 10$}\\ \cline{4-7} 
		      &    &   &    $\mathcal{E}_2[u]$(1) & order &$\mathcal{E}_2[u]$(10) & order \\ 
		\midrule
		$1\times 10^6$ &  &        & 0.0574 & -    & 0.0617 & -    \\ 
		$4\times 10^6$ & 0.01 & 0.01       & 0.0292& $0.49$ & 0.0316 & 0.48  \\ 
		$1.6\times 10^7$ & &     & 0.0149 & $0.49$ & 0.0163 & 0.48   \\ \midrule
		& & 0.25 & 0.0355 & - & 0.0968 & -\\
		$1\times 10^8$ & 0.01 & $0.2$ &  0.0288 & 0.94 &0.0792 & 0.90  \\ 
		&  & $0.1$   & 0.0158 & $0.88$  & 0.0422 & 0.91   \\ \midrule
		& $0.2$ &       & 0.0333 & -  & 0.0115 & -  \\ 
		
		$1\times 10^8$ & $0.15$ & 0.01      & 0.0248 & 0.99  & 0.0081 & 1.21     \\ 
		& $0.1$ &        & 0.0170 & $0.96$  & 0.0054 & 1.09   \\
		\hline 
	\end{tabular}
	\label{tableABaccuracy}
\end{table}

\section{Solving the 6-D Allen-Cahn equation}
\label{sec:ac}
This section is devoted into numerical experiments using SPM to integrate the 6-D Allen-Cahn equation with the local and nonlocal Laplacian operator. In order to visualize high-dimensional solutions, we adopt the following 1-D and 2-D projections
\begin{equation}\label{P&M}
	P(x_1,t) = \int_{\mathbb{R}^{d-1}} u(\bx,t) \D x_2 \dots \D x_d, \quad  M(x_1,x_2,t) = \int_{\mathbbm{R}^{d-2}} u(\bx,t) \D x_3 \dots \D x_d,
\end{equation}
and then use the same relative $L^2$ errors $\mathcal{E}_2[P](t)$ and $\mathcal{E}_2[M](t)$ defined in Eq.~\eqref{eq:L2} to measure the accuracy. The 1-D and 2-D projections are approximated by a uniform partition (choosing test functions $\varphi(\bx) = \mone_{Q^{\mu}}$ and $\varphi(\bx) = \mone_{Q^{\mu}\times Q^{\nu}}$),
\begin{equation}
	\begin{aligned}
		P(x_1,t) &\approx \sum_{\mu=1}^{(r_1-l_1)/h} \left(\frac{1}{N |Q^{\mu}|} \sum_{i=1}^N w_i(t) \mone_{Q^{\mu}} (\bx_i(t))\right) \mone_{Q^{\mu}} (x_1), \\
		M(x_1,x_2, t) & \approx \sum_{\mu=1}^{(r_1-l_1)/h} \sum_{\nu=1}^{(r_2-l_2)/h} \left(\frac{1}{N |Q^{\mu}| |Q^{\nu}|} \sum_{i=1}^N w_i(t) \mone_{Q^{\mu}\times Q^{\nu}} (\bx_i(t))\right) \mone_{Q^{\mu}\times Q^{\nu}} (x_1, x_2),
	\end{aligned}
\end{equation}
where $Q^{\mu} = \left[l_1+(\mu-1)h, l_1+\mu h\right]$ and $Q^{\nu} = \left[l_2+(\nu-1)h, l_2+\nu h\right]$. 

All simulations via our C++ implementations run on the High-Performance Computing Platform of Peking University: 2*Intel Xeon E5-2697A-v4 (2.60GHz, 40MB Cache, 9.6GT/s QPI Speed, 16 Cores, 32 Threads) with 256GB Memory $\times$ 16. 

\subsection{The Allen-Cahn equation}\label{sec:local ac}
\begin{figure}[htbp]
	\centering
	\begin{subfigure}[b]{0.3\textwidth}
		\centering
		\includegraphics[width=\textwidth]{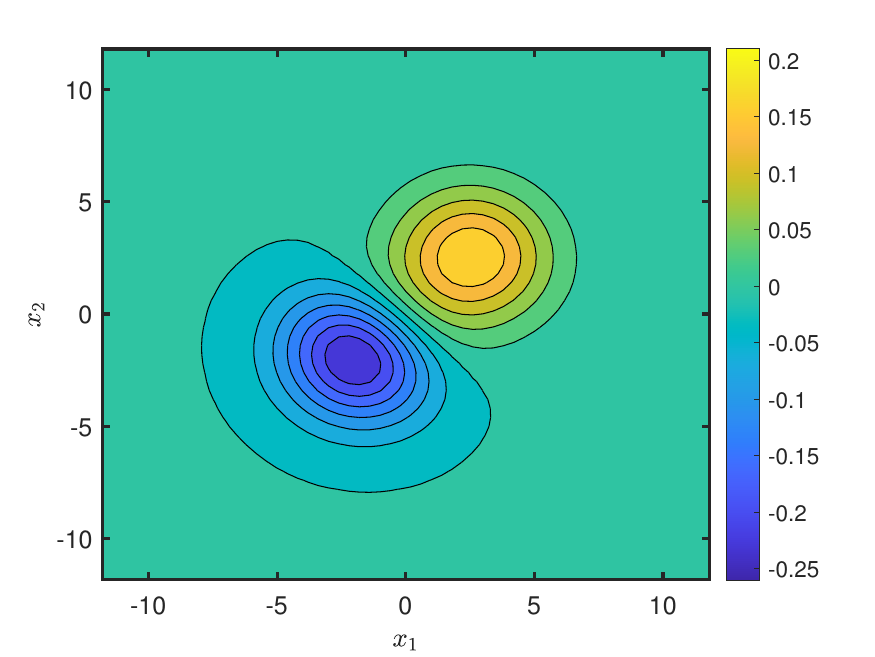}
		\caption{Reference solution.}
	\end{subfigure}
	\hfill
	\begin{subfigure}[b]{0.3\textwidth}
		\centering
		\includegraphics[width=\textwidth]{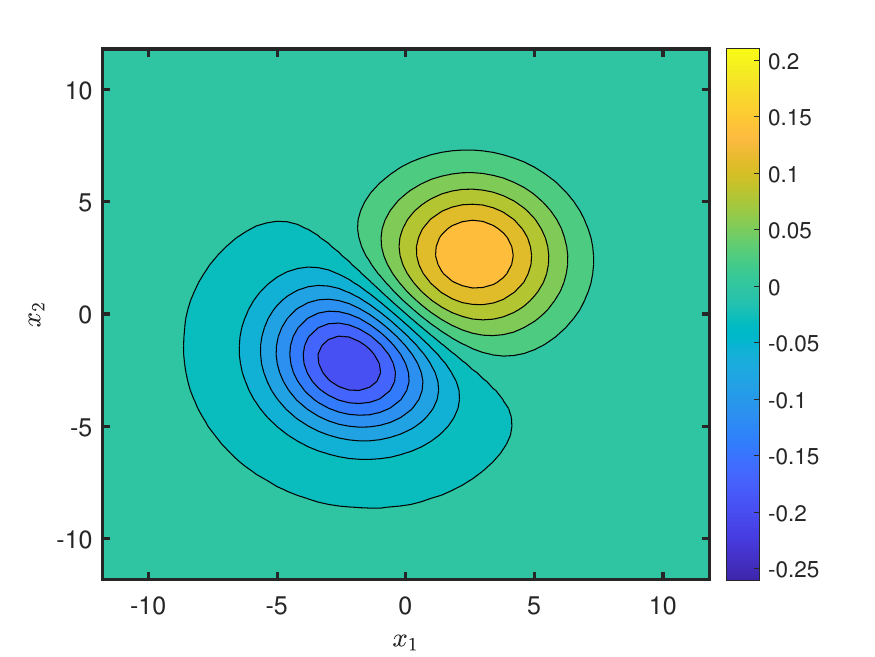}
		\caption{SPM solution with $N=4\times 10^8$.}
	\end{subfigure}
	\hfill
	\begin{subfigure}[b]{0.3\textwidth}
		\centering
		\includegraphics[width=\textwidth]{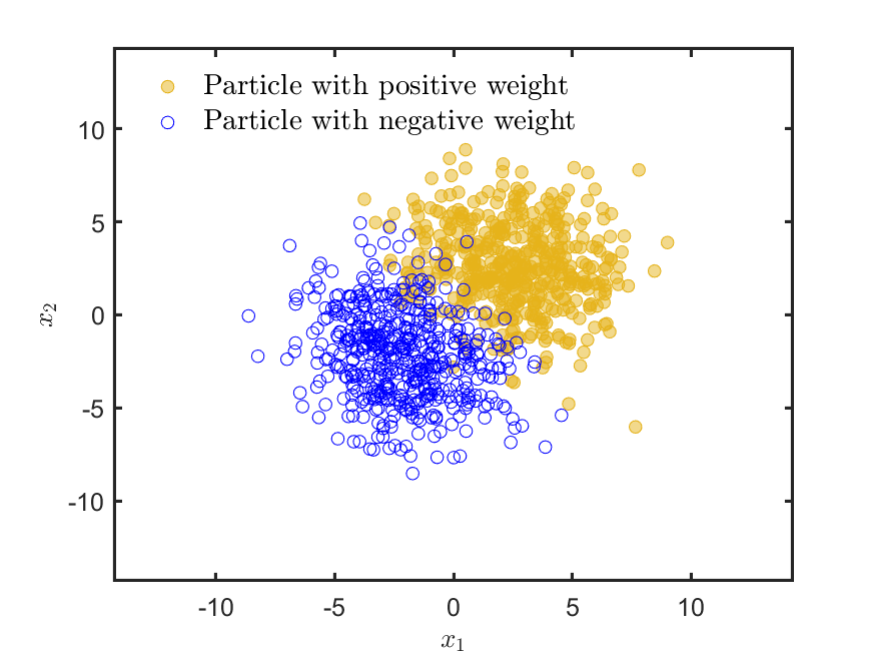}
		\caption{Particle distribution. \label{fig allencahnparticlex}}
	\end{subfigure}
	\caption{The 6-D Allen-Cahn equation with diffusion coefficient $c=1$: Filled contour plots of $M(x_1, x_2, 2)$ defined in Eq.~\eqref{P&M} for (a) the reference solution given in Eq.~\eqref{allencahnexact} and (b) the numerical solution produced by SPM with  $N=4\times 10^8$. The agreement between them is evident. SPM shows a highly adaptive characteristic since the particles with positive and negative weights concentrate in the area with positive and negative values of the solution, respectively. This is clearly demonstrated in (c) where we have randomly chosen $10^3$ particles from all samples at the final instant $T=2$ and projected their locations onto the $x_1x_2$-plane.}
	\label{fig particle following solution allencahn}
\end{figure}

\begin{figure}[htbp]
	\centering
	\begin{subfigure}[b]{0.3\textwidth}
		\centering
		\includegraphics[width=\textwidth]{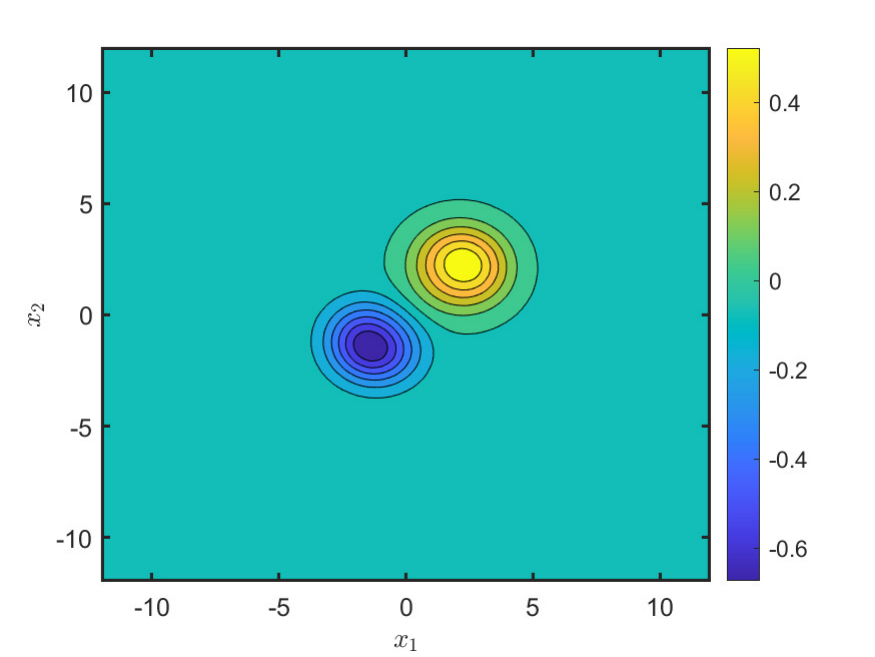}
		\caption{SPM solution for $c=0.1$.}
	\end{subfigure}
	\hfill
	\begin{subfigure}[b]{0.3\textwidth}
		\centering
		\includegraphics[width=\textwidth]{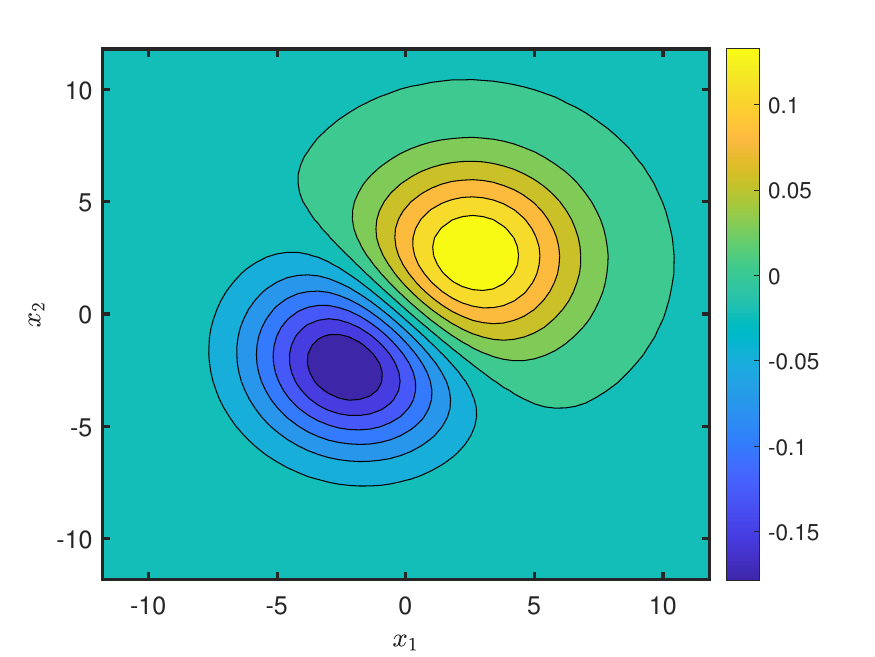}
		\caption{SPM solution for $c=1.5$.}
	\end{subfigure}
	\hfill
	\begin{subfigure}[b]{0.3\textwidth}
		\centering
		\includegraphics[width=\textwidth]{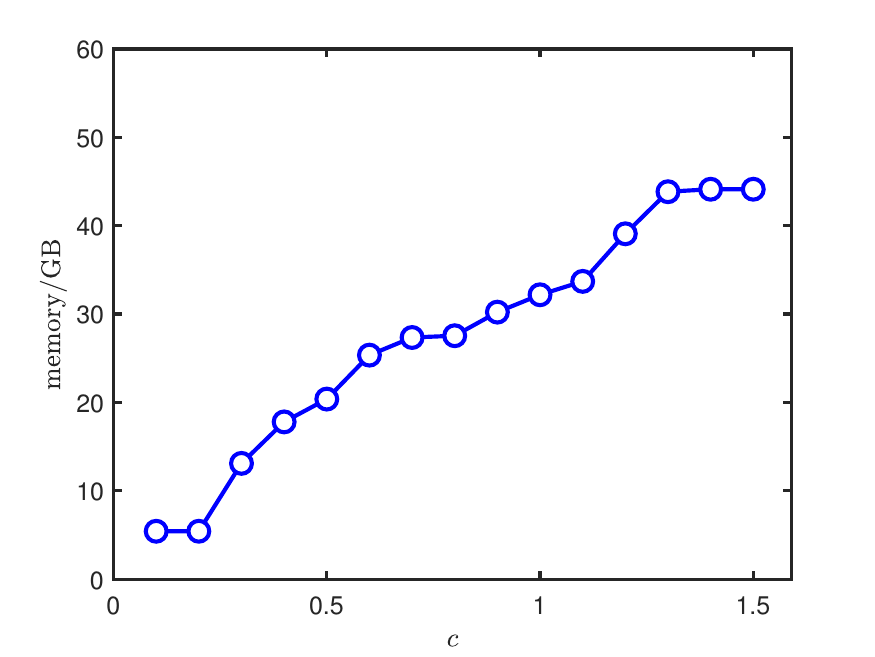}
		\caption{Memory usage against $c$.}
	\end{subfigure}
		\caption{The 6-D Allen-Cahn equation: Filled contour plots of $M(x_1, x_2, 2)$ for (a) $c=0.1$ and (b) $c=1.5$. As the diffusion coefficient $c$ increases, the support of the solution expands and the memory usage correspondingly increases as shown in (c). Here we set $N=4\times 10^8$.}
	\label{allencahn mem}
\end{figure}

We first consider the 6-D Allen-Cahn equation with local Laplacian operator

\begin{equation}\label{ac}
\frac{\partial}{\partial t}u(\vec{x},t) = c\Delta u(\vec{x},t) +  u(\vec{x},t)-{ u^3(\vec{x},t)} + {r(\vec{x},t)},\quad \bx \in \mathbb{R}^6,
\end{equation}
which allows the following analytical solution 
\begin{equation}\label{allencahnexact}
u^{\text{ref}}(\vec{x},t) = \frac{x_1+x_2}{\left(\pi (1+4ct)\right)^{d/2}}\left(\exp\left(-\frac{||\vec{x}-\vec{p}_1||_2^2}{1+4ct}\right) +2.0 \exp\left(-\frac{||\vec{x}-\vec{p}_2||_2^2}{1+4ct}\right)\right),
\end{equation}
where $\vec{p}_1 = \left(2,2,0,\dots,0\right)$, $\vec{p}_2 = \left(-1,-1,0,\dots,0\right)$, and $c$ gives the diffusion coefficient which is set to be different values in experiments. The initial data $u_0(\bx) = u^{\text{ref}}(\vec{x},0)$ and let $r(\vec{x},t)$ be the remaining term after substituting the analytical solution Eq.~\eqref{allencahnexact} into Eq.~\eqref{ac}. 
We set $\tau = 0.1$ in the Lawson-Euler scheme \eqref{Lawson_Euler} and $h = 0.4$ in the piecewise constant reconstruction via VUG adopted by Algorithm~\ref{alg SPM}.  Table~\ref{allencahnLawsonEulertable}, Figures~\ref{fig particle following solution allencahn} and \ref{allencahn mem} present the numerical results at the final time $T = 2$.

\begin{table}[htbp]
\centering
\caption{The 6-D Allen-Cahn equation with diffusion coefficient $c=1$: Total wall time with 8 cores in parallel and the relative $L^2$ errors for different sample sizes.}
\begin{tabular}{cccc}
\toprule
$N$    & $\mathcal{E}_2[P](2)$& $\mathcal{E}_2[M](2)$  & Time/h    \\ \midrule
$1\times10^8$ & 0.249& 0.254 & 0.27    \\ 
$2\times10^8$ & 0.154& 0.159 & 0.53    \\ 
$4\times10^8$ & 0.128& 0.145 & 0.96  \\ \bottomrule
\end{tabular}
\label{allencahnLawsonEulertable}
\end{table}

For the diffusion coefficient $c=1$,  
SPM takes about an hour using the distributed parallel technology via MPI with 8 cores to evolve 
the 6-D Allen-Cahn equation~\eqref{ac} until $T=2$ while maintaining the relative $L^2$  errors $\mathcal{E}_2[P](2)$ and $\mathcal{E}_2[M](2)$ less than 15\% (see Table \ref{allencahnLawsonEulertable}).
Such high efficiency fully benefits from the intrinsic adaptive characteristic of SPM. 
Figure~\ref{fig particle following solution allencahn} presents the filled contour plots of $M(x_1, x_2, 2)$ 
for the reference solution given in Eq.~\eqref{allencahnexact}, which coincides with the SPM solution produced with $N=4\times 10^8$ particles. All the particles concentrate
in the important areas in an adaptive manner as clearly demonstrated  in the last plot of Figure~\ref{fig particle following solution allencahn}, where we have randomly chosen $10^3$ particles at the final time and projected their locations into the $x_1x_2$-plane. In the area where the solution vanishes, there are always few particles. 
In addition, as displayed in Figure~\ref{allencahn mem}, when the diffusion coefficient $c$ increases, the support of the solution expands  and thus the memory usage correspondingly increases.  For instance,  
the last plot of Figure~\ref{allencahn mem} shows that the memory consumption increases from $5.4$ GB to $44.1$ GB as 
$c$ increases from $0.1$ to $1.5$.

\subsection{The Allen-Cahn equation with the nonlocal Laplacian operator}\label{sec fraction ac}

\begin{figure}[htbp]
	\centering
	\begin{subfigure}[b]{0.49\textwidth}
		\centering
		\includegraphics[width=\textwidth]{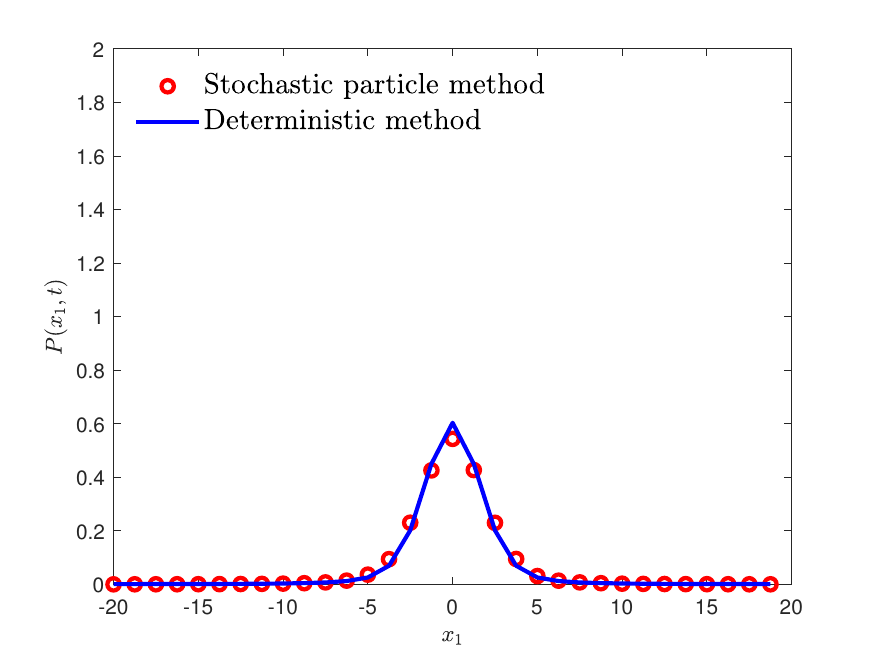}
		\caption{$t = 1, h = 0.2$}
	\end{subfigure}
	\begin{subfigure}[b]{0.49\textwidth}
		\centering
		\includegraphics[width=\textwidth]{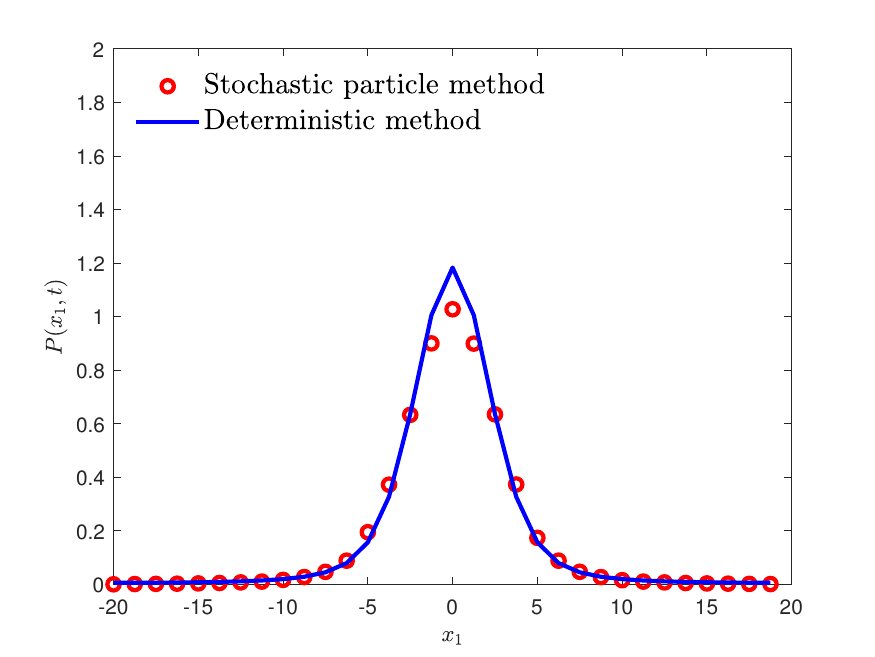}
		\caption{$t = 2, h = 0.2$}
	\end{subfigure}
	\begin{subfigure}[b]{0.49\textwidth}
		\centering
		\includegraphics[width=\textwidth]{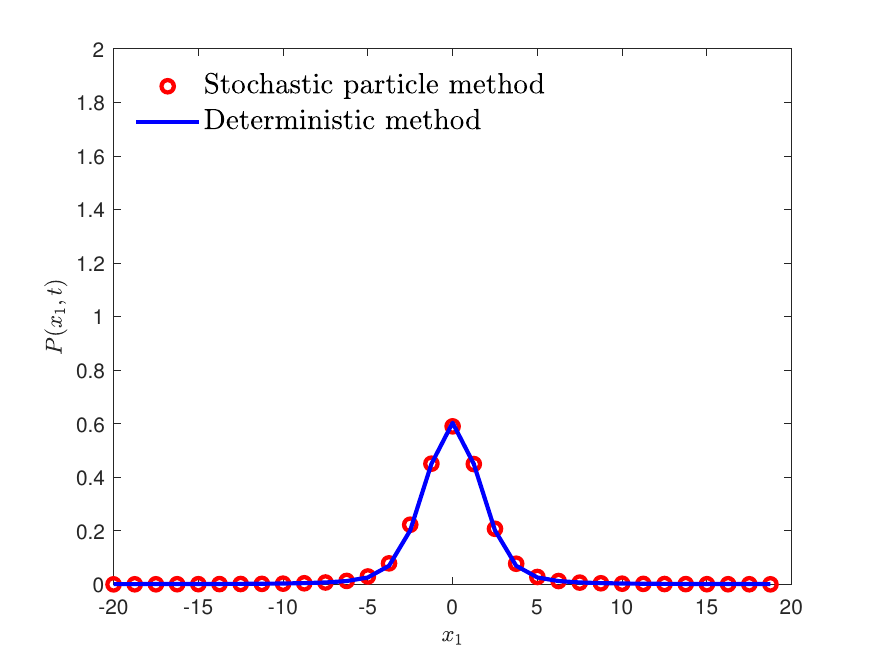}
		\caption{$t = 1, h = 0.1$}
	\end{subfigure}
	\begin{subfigure}[b]{0.49\textwidth}
		\centering
		\includegraphics[width=\textwidth]{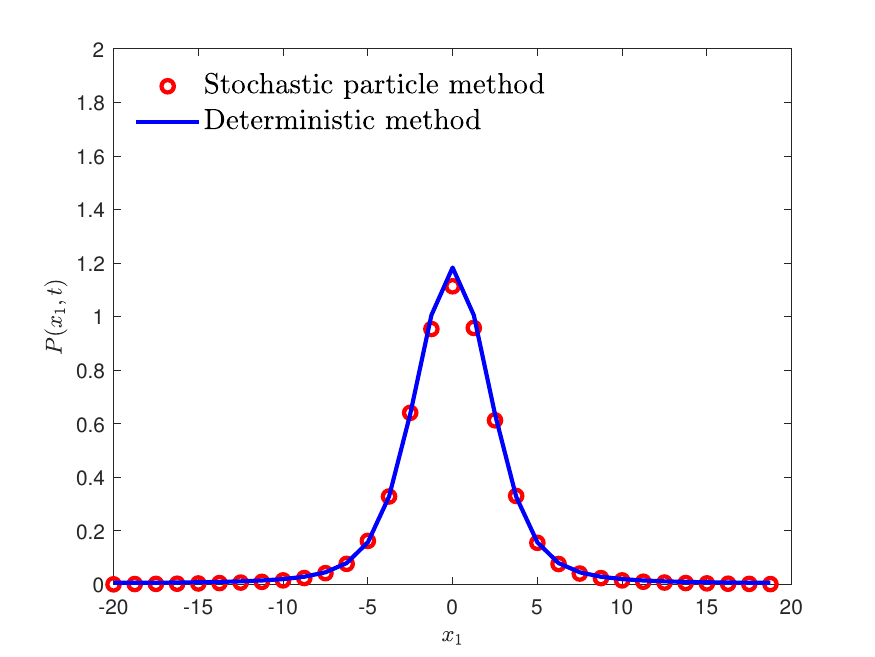}
		\caption{$t = 2, h = 0.1$}
	\end{subfigure}
	\caption{The 6-D Allen-Cahn equation with the nonlocal Laplacian operator: Plots of $P(x_1,t)$ defined in Eq.~\eqref{P&M} at $t=1, 2$ for the SPM solutions with $N = 1\times 10^8$ and $\tau = 0.01$ against the reference ones. In (c) and (d), the SPM solution with $h = 0.1$ match the reference solution better than the SPM solution with $h = 0.2$ in (a) and (b), particularly at the peak.}
	\label{fig 6d frac}
\end{figure}

\begin{figure}[htbp]
	\centering
	\begin{subfigure}[b]{0.45\textwidth}
		\centering
		\includegraphics[width=\textwidth]{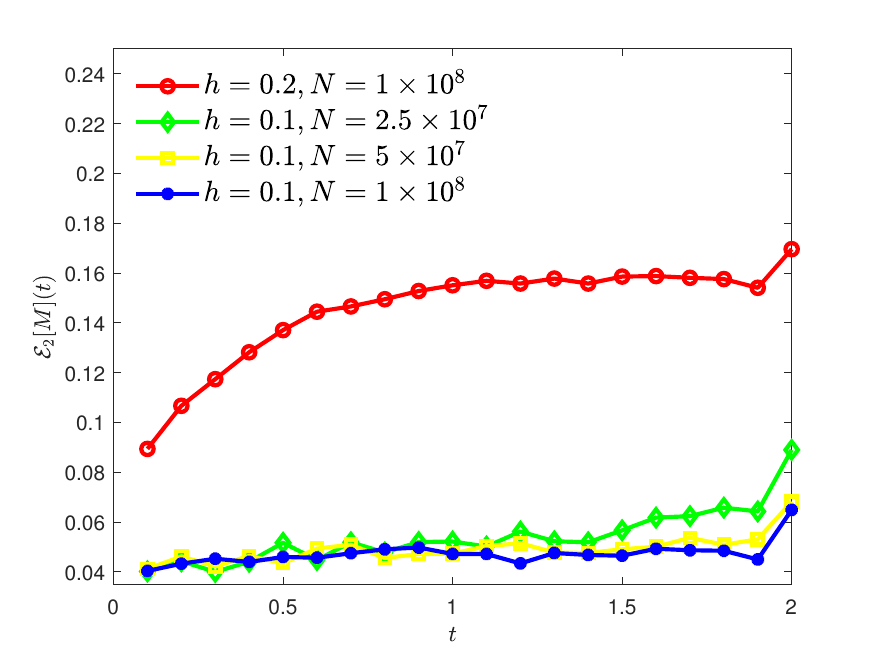}
	\end{subfigure}
	\caption{The 6-D Allen-Cahn equation with the nonlocal Laplacian operator: The time evolution of the relative $L^2$ errors $\mathcal{E}_2[M](t)$ under different $h$ and $N$. The errors with $h = 0.2$ are significantly larger than that with $h = 0.1$, which may be due to the insufficient resolution at the peak as already shown in Figure \ref{fig 6d frac}. The errors with $h = 0.1$ generally decrease as $N$ increases. }
	\label{fig 6d frac error}
\end{figure}

Next, we consider the Allen-Cahn equation with the nonlocal Laplacian operator
\begin{equation}\label{6d nonlocal ac}
	\begin{aligned}
		\frac{\partial}{\partial t} u(\boldsymbol{x}, t)&=-(-\Delta)^{\alpha / 2} u(\boldsymbol{x}, t)+u(\vec{x},t)-u^3(\vec{x},t),\quad \bx \in \mathbbm{R}^6,\\
		u(\bx,0) &= \frac{1}{(2\pi)^{3}}\exp(-\Vert \bx \Vert^2/2),
	\end{aligned}
\end{equation}
and set $\alpha = 1.5$. The random walk algorithm for the nonlocal operator given in Algorithm~\ref{alg fractional random walk} runs with the cut-off parameter $\epsilon=0.005$. To measure the accuracy, a deterministic solver attempts to solve the 6-D problem, which uses the second-order operator splitting \cite{hansen2012second} and spectral approximations are employed to produce the reference solution on a $32^6$ uniform grid.
Figure \ref{fig 6d frac} plots $P(x_1,1)$ and $P(x_1,2)$ with $h = 0.1, 0.2$ and clearly shows that a smaller $h$ enables the SPM solutions to fit better with the reference ones at the peak. This is also confirmed in Figure \ref{fig 6d frac error} which displays the time evolution of $\mathcal{E}_2[M](t)$ under different $h$ and $N$. As expected, it further shows there that the increase of sample size $N$ may generally decrease the errors.  Table \ref{table 6d nonlinear} presents the relative $L^2$ errors and total wall time for different particle number $N$. It takes nearly an hour to make $\mathcal{E}_2[P](2)$ and $\mathcal{E}_2[M](2)$ achieve 4.92\% and 6.49\% respectively.




\begin{table}[H]
	\centering
	\caption{The 6-D Allen-Cahn equation with the nonlocal Laplacian operator: Total wall time with 32 cores in parallel and the relative $L^2$ errors for different sample sizes. Here we set $\tau = 0.01, h = 0.1$.}
	\begin{tabular}{cccc}
		\toprule
		$N $       & $\mathcal{E}_2[P](2)$  & $\mathcal{E}_2[M](2)$ & Time/h   \\ \midrule
		$2.5\times 10^7$ & 0.0746&0.0890 & 0.28  \\ 
		 $5\times 10^7$ & 0.0528&0.0682 & 0.66\\ 
		 $1\times 10^8$ & 0.0492&0.0649 & 1.05   \\ \bottomrule
	\end{tabular}
	\label{table 6d nonlinear}
\end{table}

\section{Solving the 7-D Hamiltonian-Jacobi-Bellman equation}
\label{sec:hjb}

\begin{figure}[H]
	\centering
	\begin{subfigure}[b]{0.45\textwidth}
		\centering
		\includegraphics[width=\textwidth]{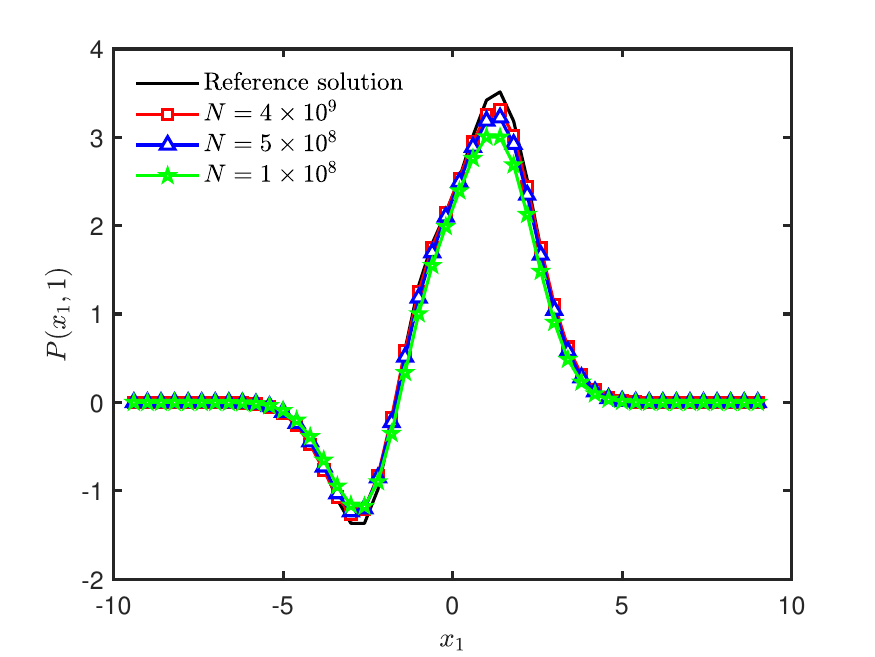}
		\caption{$t=1$.}
	\end{subfigure}
	\hfill
	\begin{subfigure}[b]{0.45\textwidth}
		\centering
		\includegraphics[width=\textwidth]{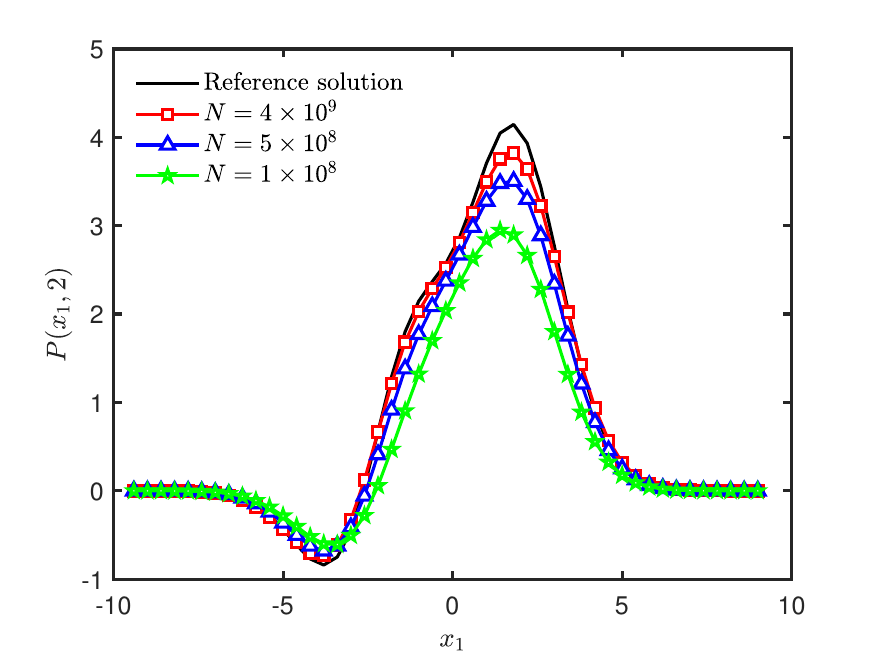}
		\caption{$t=2$.}
		\label{hjbt2}
	\end{subfigure}
	\caption{The 7-D HJB equation: Plots of $P(x_1,t)$ defined in Eq.~\eqref{P&M} at $t=1, 2$ for the SPM solutions with $N=1\times 10^8, 5\times 10^8$ and $4\times 10^9$ against the reference solution given in Eq.~\eqref{HJBexact}.}
	\label{fig HJB1dproj}
\end{figure}

\begin{figure}[H]
\centering
\begin{subfigure}[b]{0.49\textwidth}
\centering
\includegraphics[width=\textwidth]{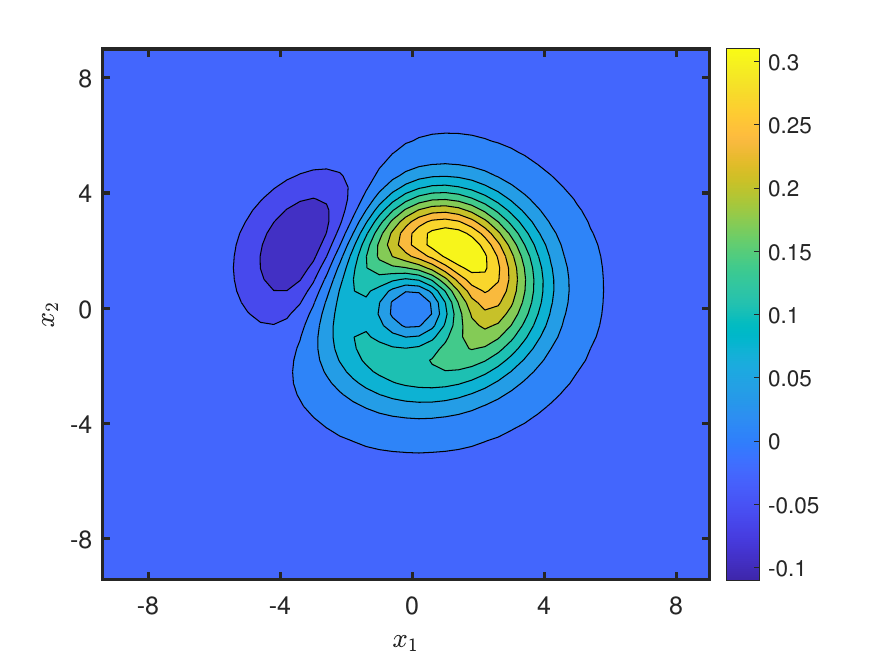}
\caption{Reference solution.}
\end{subfigure}
\hfill
\begin{subfigure}[b]{0.49\textwidth}
\centering
\includegraphics[width=\textwidth]{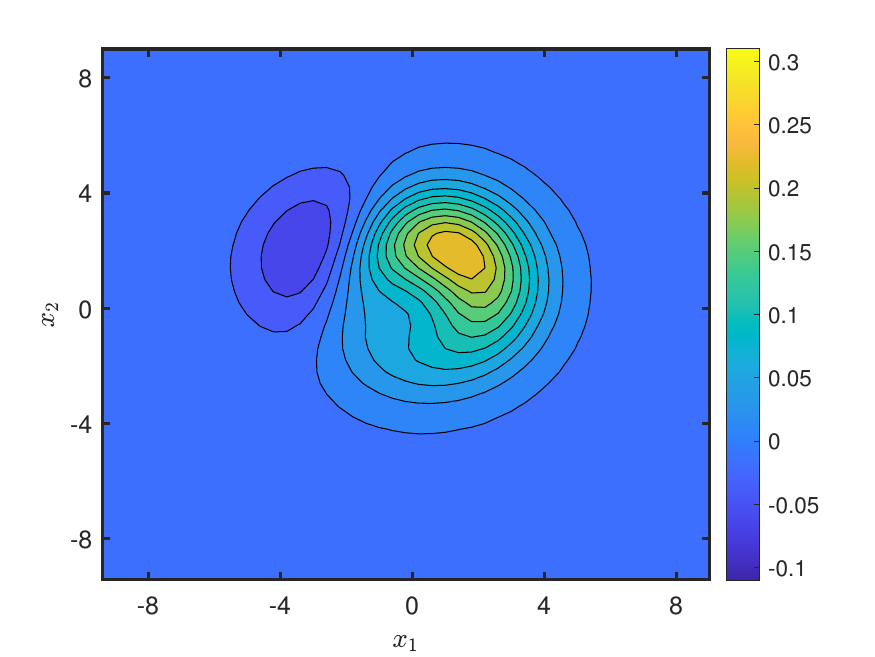}
\caption{$N=1\times 10^8$.}
\end{subfigure}
\hfill
\begin{subfigure}[b]{0.49\textwidth}
\centering
\includegraphics[width=\textwidth]{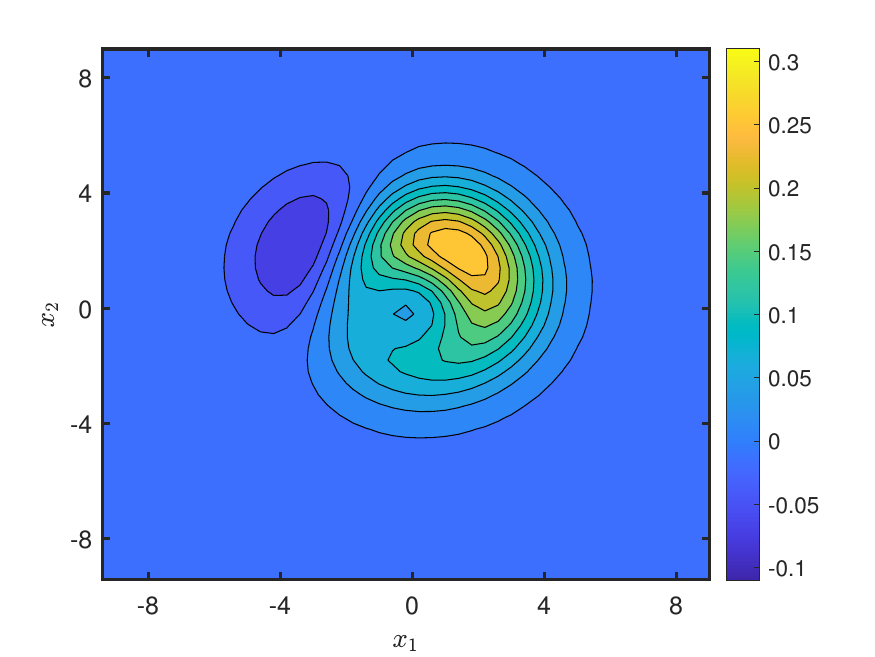}
\caption{$N=5\times 10^8$.}
\end{subfigure}
\hfill
\begin{subfigure}[b]{0.49\textwidth}
\centering
\includegraphics[width=\textwidth]{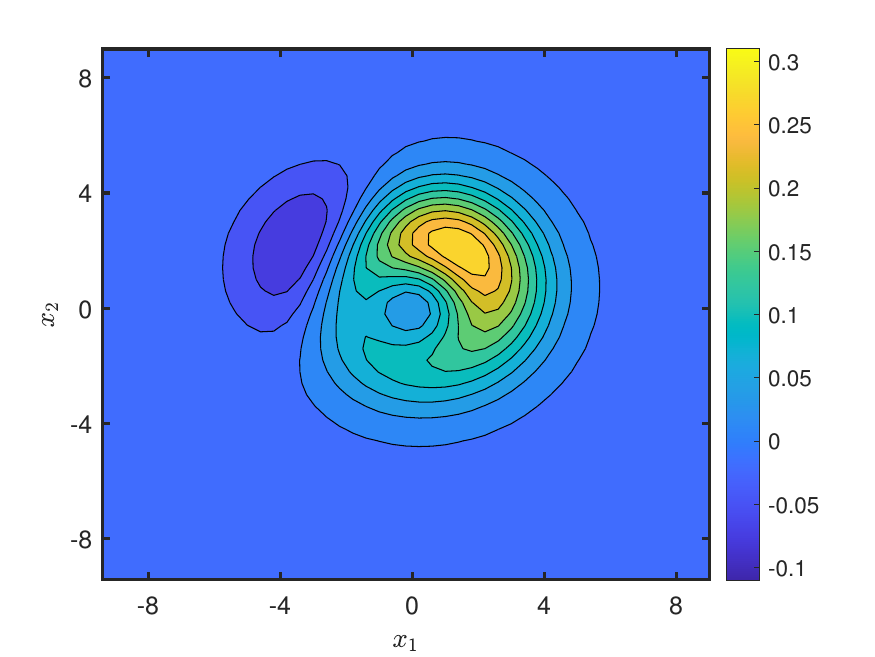}
\caption{$N=4\times 10^9$.}
\end{subfigure}
\caption{The 7-D HJB equation: Filled contour plots of $M(x_1, x_2, 2)$ defined in Eq.~\eqref{P&M} for (a) the reference solution given in Eq.~\eqref{HJBexact} and the numerical solutions produced by SPM with (b) $N=1\times 10^8$, (c) $5\times 10^8$ and  (d) $4\times 10^9$.}
\label{HJBLawsonEulerfig}
\end{figure}

Finally, we apply SPM to integrate the 7-D HJB equation
\begin{equation}
\frac{\partial }{\partial t}u(\vec{x},t) = c\Delta u(\vec{x},t) + ||\vec{\nabla }u(\vec{x},t)||^2_2 + {r(\vec{x},t)}, \quad \bx \in \mathbb{R}^7,
\end{equation}
which admits the following analytical solution 
\begin{equation}\label{HJBexact}
u^\text{ref}(\vec{x},t) = \frac{x_1^2+x_2^2}{\left(\pi (1+4ct)\right)^{d/2}}\left(2.5\exp\left(-\frac{||\vec{x}-\vec{p}_1||_2^2}{1+4ct}\right) - \exp\left(-\frac{||\vec{x}-\vec{p}_2||_2^2}{1+4ct}\right)\right),
\end{equation}
where $\vec{p}_1 = \left(0,0.6,0,\dots,0\right)$, $\vec{p}_2 = \left(-1,1,0,\dots,0\right)$ and $c=0.5$. We adopt the same $T, \tau, h$ as in Section \ref{sec:local ac}. 
Figures~\ref{fig HJB1dproj} and \ref{HJBLawsonEulerfig} plot $P(x_1,t)$ and $M(x_1, x_2, 2)$, respectively, defined in Eq.~\eqref{P&M} for the numerical solutions produced by SPM with $N=1\times 10^8, 5\times 10^8$ and $4\times 10^9$ against the reference solution given in Eq.~\eqref{HJBexact}. Table~\ref{HJBLawsonEulertable} shows the corresponding relative $L^2$ errors, memory usage as well as the total wall time.
SPM can systematically increase the accuracy. We are able to observe there that, when the number of particles $N$ increases from $1\times 10^8$ to $4\times 10^9$, the error $\mathcal{E}_2[P](2)$ decreases from 30.8\% to 6.1\% and $\mathcal{E}_2[M](2)$ from 34.5\% to 10.5\%,
whereas the computational time with 8 cores shows a sharp increase from about half an hour to $10$ hours 
and so does the the memory usage.


%

\begin{table}[H]
	\centering
	\caption{The 7-D HJB equation: Total wall time with 8 cores in parallel, memory usage and the relative $L^2$ errors under different sample sizes. Here $P_{\text{grad}}(x_1,t) = \int_{\mathbb{R}^{d-1}} \Vert \vec{\nabla} u(\bx,t)\Vert_2^2 \,\,\D x_2\dots \D x_d$ and $M_{\text{grad}}(x_1,x_2,t) = \int_{\mathbb{R}^{d-2}} \Vert \vec{\nabla} u(\bx,t)\Vert_2^2 \,\,\D x_3\dots \D x_d$}.
	\begin{tabular}{cccc}
		\toprule
		$N$ & $1\times 10^8$ & $5\times 10^8$ & $4\times 10^9$ \\ \midrule
		$\mathcal{E}_2[P](2)$ & 0.308 & 0.145 & 0.061 \\
		$\mathcal{E}_2[M](2)$ & 0.345 & 0.184 & 0.105 \\
		$\mathcal{E}_2[P_{\text{grad}}](1)$ & 0.166 & 0.153 & 0.107 \\
		$\mathcal{E}_2[M_{\text{grad}}](1)$ & 0.246 & 0.216 & 0.188 \\
		$\mathcal{E}_2[P_{\text{grad}}](2)$ & 0.330 & 0.235 & 0.170 \\
		$\mathcal{E}_2[M_{\text{grad}}](2)$ & 0.484 & 0.335 & 0.290 \\
		Memory/GB & 17.43 & 73.58 & 200.07 \\
		Time/h & 0.39 & 1.60 & 9.81 \\
		\bottomrule
	\end{tabular}
\label{HJBLawsonEulertable}
\end{table}

\section{Conclusion and discussion}
\label{sec:con}


This work proposes a stochastic particle method (SPM) to moderately high-dimensional nonlinear PDEs. Starting from the weak formulation of the Lawson-Euler scheme,  SPM uses  real-valued weighted particles to approximate the high-dimensional solution in the weak sense,
and updates the particle locations and weights at each time step.
In the spirit of importance sampling, the relocating technique is proposed to make particle locations match the shape of the solution and reduce the stochastic variances. A piecewise constant reconstruction with virtual uniform grid (VUG),  which makes full use of the 
adaptive characteristics of stochastic particles moving with the solution, 
is adopted for obtaining the values of nonlinear terms. Numerical experiments in solving the 6-D Allen-Cahn equation and the 7-D Hamiltonian-Jacobi-Bellman equation demonstrate the potential of SPM in solving moderately high-dimensional nonlinear PDEs efficiently while maintaining a reasonable accuracy.

Apart from SPM, there have been fruitful developments in solving high-dimensional problems recently, such as deep learning \cite{han2017deep,nusken2021interpolating,GaoWang2023}, tensor trains \cite{bachmayr2023low,bachmayr2016tensor,dolgov2021tensor, tang2024solving, richter2024continuous} and sparse grids \cite{Griebel2004Sparsegrids, kormann2016sparse}. 
It definitely contributes to a better understanding of SPM to undertake certain and preliminary discussions regarding the potential advantages and disadvantages among all these bold attempts. 
\begin{itemize}
	\item Compared with deep learning: There is no learning and hyper-parameter tuning in SPM. The impact of parameters on errors is crystal clear, and errors can be systematically reduced. But currently SPM only deals with PDEs in moderately large dimension, while deep learning might potentially deal with problems in much higher dimension.

	\item Compared with tensor trains: SPM does not requires unknown solutions to have a low-rank structure. Instead, it depends on the intrinsic adaptive features of particles to economize the computational cost on the precondition of maintaining the accuracy. When the low-rank structure is present, the tensor trains method may be more efficient \cite{dolgov2021tensor}.

	\item Compared with sparse grids: SPM automatically avoids unnecessary storage through the inherent adaptability of stochastic particles. By contrast, in sparse grids, the storage is reduced by  omitting coefficients whose absolute values are smaller than a given tolerance \cite{Griebel2004Sparsegrids}. 
\end{itemize}


\par Solving certain linear problems by SPM might be immune to the CoD (see the high-dimensional linear example in \ref{appendix A}). Reconstructing the solution values from weighted particles, which is a relatively independent topic, is of high complexity in SPM, and it poses a challenging problem for SPM in solving nonlinear PDEs with higher dimensions. In spite of the uniform piecewise constant reconstruction used in this work, the adaptive piecewise constant reconstruction, neural networks and radial basis functions may potentially be combined with SPM, which will be investigated in the future work.

Finally, we would like to point out that the proposed SPM provides a particle-based framework for general nonlinear PDEs without transparent probabilistic interpretation. Actually, the variational form of any discretization scheme,  including finite difference, finite element, spectral methods, can be cast into a stochastic representation.  In this sense, SPM inherits the celebrated theory from the classical numerical schemes that have been well-established over the last 100 years,  and extends their boundaries. In other words, SPM grows on sound and solid mathematical ground.

\section*{Acknowledgement}
This research was supported by the National Natural Science Foundation of China (Nos.~12325112, 12101018, 12288101), 
the Fundamental Research Funds for the Central Universities (No.~310421125) and the High-performance Computing Platform of Peking University. The authors are sincerely grateful to Haoyang Liu and Shuyi Zhang at Peking University for their technical support on the computing environment.



\begin{thebibliography}{54}
	\expandafter\ifx\csname natexlab\endcsname\relax\def\natexlab#1{#1}\fi
	\providecommand{\url}[1]{\texttt{#1}}
	\providecommand{\href}[2]{#2}
	\providecommand{\path}[1]{#1}
	\providecommand{\DOIprefix}{doi:}
	\providecommand{\ArXivprefix}{arXiv:}
	\providecommand{\URLprefix}{URL: }
	\providecommand{\Pubmedprefix}{pmid:}
	\providecommand{\doi}[1]{\href{http://dx.doi.org/#1}{\path{#1}}}
	\providecommand{\Pubmed}[1]{\href{pmid:#1}{\path{#1}}}
	\providecommand{\bibinfo}[2]{#2}
	\ifx\xfnm\relax \def\xfnm[#1]{\unskip,\space#1}\fi
	\bibitem[{Smyth et~al.(1998)Smyth, Parker, and Taylor}]{smyth1998numerical}
	\bibinfo{author}{E.~S. Smyth}, \bibinfo{author}{J.~S. Parker},
	\bibinfo{author}{K.~T. Taylor},
	\newblock \bibinfo{title}{Numerical integration of the time-dependent
		{S}chr{\"o}dinger equation for laser-driven {H}elium},
	\newblock \bibinfo{journal}{Commun. Comput. Phys.} \bibinfo{volume}{114}
	(\bibinfo{year}{1998}) \bibinfo{pages}{1--14}.
	\bibitem[{Xiong and Shao(2024)}]{XiongShao2020Overcoming}
	\bibinfo{author}{Y.~Xiong}, \bibinfo{author}{S.~Shao},
	\newblock \bibinfo{title}{Overcoming the numerical sign problem in the {W}igner
		dynamics via adaptive particle annihilation},
	\newblock \bibinfo{journal}{SIAM J. Sci. Comput.} \bibinfo{volume}{46}
	(\bibinfo{year}{2024}) \bibinfo{pages}{B107--B136}.
	\bibitem[{Kormann et~al.(2019)Kormann, Reuter, and
		Rampp}]{kormann2019massively}
	\bibinfo{author}{K.~Kormann}, \bibinfo{author}{K.~Reuter},
	\bibinfo{author}{M.~Rampp},
	\newblock \bibinfo{title}{A massively parallel semi-{L}agrangian solver for the
		six-dimensional {V}lasov-{P}oisson equation},
	\newblock \bibinfo{journal}{Int. J. High Perform. Comput. Appl.}
	\bibinfo{volume}{33} (\bibinfo{year}{2019}) \bibinfo{pages}{924--947}.
	\bibitem[{Dimarco et~al.(2018)Dimarco, Loub{\`e}re, Narski, and
		Rey}]{dimarco2018efficient}
	\bibinfo{author}{G.~Dimarco}, \bibinfo{author}{R.~Loub{\`e}re},
	\bibinfo{author}{J.~Narski}, \bibinfo{author}{T.~Rey},
	\newblock \bibinfo{title}{An efficient numerical method for solving the
		{B}oltzmann equation in multidimensions},
	\newblock \bibinfo{journal}{J. Comput. Phys.} \bibinfo{volume}{353}
	(\bibinfo{year}{2018}) \bibinfo{pages}{46--81}.
	\bibitem[{Bellman(1957)}]{bk:Bellman1957}
	\bibinfo{author}{R.~E. Bellman}, \bibinfo{title}{Dynamic Programming},
	\bibinfo{publisher}{Princeton University Press}, \bibinfo{year}{1957}.
	\bibitem[{Xiong et~al.(2023)Xiong, Zhang, and Shao}]{xiong2022characteristic}
	\bibinfo{author}{Y.~Xiong}, \bibinfo{author}{Y.~Zhang},
	\bibinfo{author}{S.~Shao},
	\newblock \bibinfo{title}{A characteristic-spectral-mixed scheme for
		six-dimensional {W}igner-{C}oulomb dynamics},
	\newblock \bibinfo{journal}{SIAM J. Sci. Comput.} \bibinfo{volume}{45}
	(\bibinfo{year}{2023}) \bibinfo{pages}{B906--B931}.
	\bibitem[{Jiang et~al.(2017)Jiang, Zhang, and Shi}]{jiang2017stability}
	\bibinfo{author}{K.~Jiang}, \bibinfo{author}{P.~Zhang},
	\bibinfo{author}{A.~Shi},
	\newblock \bibinfo{title}{Stability of icosahedral quasicrystals in a simple
		model with two-length scales},
	\newblock \bibinfo{journal}{J. Phys.: Condens. Matter} \bibinfo{volume}{29}
	(\bibinfo{year}{2017}) \bibinfo{pages}{124003}.
	\bibitem[{Bachmayr(2023)}]{bachmayr2023low}
	\bibinfo{author}{M.~Bachmayr},
	\newblock \bibinfo{title}{Low-rank tensor methods for partial differential
		equations},
	\newblock \bibinfo{journal}{Acta Numer.} \bibinfo{volume}{32}
	(\bibinfo{year}{2023}) \bibinfo{pages}{1--121}.
	\bibitem[{Bachmayr et~al.(2016)Bachmayr, Schneider, and
		Uschmajew}]{bachmayr2016tensor}
	\bibinfo{author}{M.~Bachmayr}, \bibinfo{author}{R.~Schneider},
	\bibinfo{author}{A.~Uschmajew},
	\newblock \bibinfo{title}{Tensor networks and hierarchical tensors for the
		solution of high-dimensional partial differential equations},
	\newblock \bibinfo{journal}{Found. Comput. Math.} \bibinfo{volume}{16}
	(\bibinfo{year}{2016}) \bibinfo{pages}{1423--1472}.
	\bibitem[{Dolgov et~al.(2021)Dolgov, Kalise, and Kunisch}]{dolgov2021tensor}
	\bibinfo{author}{S.~Dolgov}, \bibinfo{author}{D.~Kalise},
	\bibinfo{author}{K.~K. Kunisch},
	\newblock \bibinfo{title}{Tensor decomposition methods for high-dimensional
		{H}amilton-{J}acobi-{B}ellman equations},
	\newblock \bibinfo{journal}{SIAM J. Sci. Comput.} \bibinfo{volume}{43}
	(\bibinfo{year}{2021}) \bibinfo{pages}{A1625--A1650}.
	\bibitem[{Tang and Ying(2024)}]{tang2024solving}
	\bibinfo{author}{X.~Tang}, \bibinfo{author}{L.~Ying},
	\newblock \bibinfo{title}{Solving high-dimensional {F}okker-{P}lanck equation
		with functional hierarchical tensor},
	\newblock \bibinfo{journal}{J. Comput. Phys.} \bibinfo{volume}{511}
	(\bibinfo{year}{2024}) \bibinfo{pages}{113110}.
	\bibitem[{Richter et~al.(2024)Richter, Sallandt, and
		N{\"u}sken}]{richter2024continuous}
	\bibinfo{author}{L.~Richter}, \bibinfo{author}{L.~Sallandt},
	\bibinfo{author}{N.~N{\"u}sken},
	\newblock \bibinfo{title}{From continuous-time formulations to discretization
		schemes: tensor trains and robust regression for {BSDE}s and parabolic
		{PDE}s},
	\newblock \bibinfo{journal}{J. Mach. Learn. Res.} \bibinfo{volume}{25}
	(\bibinfo{year}{2024}) \bibinfo{pages}{1--40}.
	\bibitem[{Li et~al.(2013)Li, Chen, and Chen}]{LiChenChen2013}
	\bibinfo{author}{M.~Li}, \bibinfo{author}{W.~Chen}, \bibinfo{author}{C.~S.
		Chen},
	\newblock \bibinfo{title}{The localized {RBF}s collocation methods for solving
		high dimensional {PDEs}},
	\newblock \bibinfo{journal}{Eng. Anal. Bound. Elem.} \bibinfo{volume}{37}
	(\bibinfo{year}{2013}) \bibinfo{pages}{1300--1304}.
	\bibitem[{Giles(2008)}]{Giles2008}
	\bibinfo{author}{M.~B. Giles},
	\newblock \bibinfo{title}{Multilevel {M}onte {C}arlo path simulation},
	\newblock \bibinfo{journal}{Oper. Res.} \bibinfo{volume}{56}
	(\bibinfo{year}{2008}) \bibinfo{pages}{607--617}.
	\bibitem[{Wang and Sloan(2005)}]{WangSloan2005}
	\bibinfo{author}{X.~Wang}, \bibinfo{author}{I.~H. Sloan},
	\newblock \bibinfo{title}{Why are high-dimensional finance problems often of
		low effective dimension?},
	\newblock \bibinfo{journal}{SIAM J. Sci. Comput.} \bibinfo{volume}{27}
	(\bibinfo{year}{2005}) \bibinfo{pages}{159--183}.
	\bibitem[{Yan et~al.(2013)Yan, Cai, and Zeng}]{YanCaiZeng2013}
	\bibinfo{author}{C.~Yan}, \bibinfo{author}{W.~Cai}, \bibinfo{author}{X.~Zeng},
	\newblock \bibinfo{title}{A parallel method for solving {L}aplace equations
		with {D}irichlet data using local boundary integral equations and random
		walks},
	\newblock \bibinfo{journal}{SIAM J. Sci. Comput.} \bibinfo{volume}{35}
	(\bibinfo{year}{2013}) \bibinfo{pages}{B868--B889}.
	\bibitem[{Su et~al.(2023)Su, Xu, and Sheng}]{su2023new}
	\bibinfo{author}{B.~Su}, \bibinfo{author}{C.~Xu}, \bibinfo{author}{C.~Sheng},
	\newblock \bibinfo{title}{A new ‘walk on spheres’ type method for
		fractional diffusion equation in high dimensions based on the {F}eynman-{K}ac
		formulas},
	\newblock \bibinfo{journal}{Appl. Math. Lett.} \bibinfo{volume}{141}
	(\bibinfo{year}{2023}) \bibinfo{pages}{108597}.
	\bibitem[{Bayer et~al.(2023)Bayer, Eigel, Sallandt, and
		Trunschke}]{BayerEigelSallandtTrunschke2023}
	\bibinfo{author}{C.~Bayer}, \bibinfo{author}{M.~Eigel},
	\bibinfo{author}{L.~Sallandt}, \bibinfo{author}{P.~Trunschke},
	\newblock \bibinfo{title}{Pricing high-dimensional {B}ermudan options with
		hierarchical tensor formats},
	\newblock \bibinfo{journal}{SIAM J. Financ. Math.} \bibinfo{volume}{14}
	(\bibinfo{year}{2023}) \bibinfo{pages}{383--406}.
	\bibitem[{E et~al.(2020)E, Ma, and Wu}]{EMaWu2019}
	\bibinfo{author}{W.~E}, \bibinfo{author}{C.~Ma}, \bibinfo{author}{L.~Wu},
	\newblock \bibinfo{title}{Machine learning from a continuous viewpoint},
	\newblock \bibinfo{journal}{Sci. China Math.} \bibinfo{volume}{63}
	(\bibinfo{year}{2020}) \bibinfo{pages}{2233--2266}.
	\bibitem[{Han et~al.(2018)Han, Jentzen, and E}]{HanJentzenE2018}
	\bibinfo{author}{J.~Han}, \bibinfo{author}{A.~Jentzen}, \bibinfo{author}{W.~E},
	\newblock \bibinfo{title}{Solving high-dimensional partial differential
		equations using deep learning},
	\newblock \bibinfo{journal}{P. Natl. Acad. Sci. USA} \bibinfo{volume}{115}
	(\bibinfo{year}{2018}) \bibinfo{pages}{8505--8510}.
	\bibitem[{Raissi et~al.(2019)Raissi, Perdikaris, and
		Karniadakis}]{RaissiPerdikarisKarniadakis2019}
	\bibinfo{author}{M.~Raissi}, \bibinfo{author}{P.~Perdikaris},
	\bibinfo{author}{G.~E. Karniadakis},
	\newblock \bibinfo{title}{Physics-informed neural networks: A deep learning
		framework for solving forward and inverse problems involving nonlinear
		partial differential equations},
	\newblock \bibinfo{journal}{J. Comput. Phys.} \bibinfo{volume}{378}
	(\bibinfo{year}{2019}) \bibinfo{pages}{686--707}.
	\bibitem[{Hur{\'e} et~al.(2020)Hur{\'e}, Pham, and Warin}]{HurePhamWarin2020}
	\bibinfo{author}{C.~Hur{\'e}}, \bibinfo{author}{H.~Pham},
	\bibinfo{author}{X.~Warin},
	\newblock \bibinfo{title}{Deep backward schemes for high-dimensional nonlinear
		{PDEs}},
	\newblock \bibinfo{journal}{Math. Comp.} \bibinfo{volume}{89}
	(\bibinfo{year}{2020}) \bibinfo{pages}{1547--1579}.
	\bibitem[{Kast and Hesthaven(2024)}]{kast2024positional}
	\bibinfo{author}{M.~Kast}, \bibinfo{author}{J.~S. Hesthaven},
	\newblock \bibinfo{title}{Positional embeddings for solving {PDE}s with
		evolutional deep neural networks},
	\newblock \bibinfo{journal}{J. Comput. Phys.} \bibinfo{volume}{508}
	(\bibinfo{year}{2024}) \bibinfo{pages}{112986}.
	\bibitem[{Richter et~al.(2021)Richter, Sallandt, and
		N{\"u}sken}]{richter2021solving}
	\bibinfo{author}{L.~Richter}, \bibinfo{author}{L.~Sallandt},
	\bibinfo{author}{N.~N{\"u}sken},
	\newblock \bibinfo{title}{Solving high-dimensional parabolic {PDE}s using the
		tensor train format},
	\newblock in: \bibinfo{booktitle}{Proc. Int. Conf. Mach. Learn.},
	\bibinfo{year}{2021}, pp. \bibinfo{pages}{8998--9009}.
	\bibitem[{Courant et~al.(1928)Courant, Friedrichs, and Lewy}]{CFL1928}
	\bibinfo{author}{R.~Courant}, \bibinfo{author}{K.~Friedrichs},
	\bibinfo{author}{H.~Lewy},
	\newblock \bibinfo{title}{{\"U}ber die partiellen differenzengleichungen der
		mathematischen physik},
	\newblock \bibinfo{journal}{Math. Ann.} \bibinfo{volume}{100}
	(\bibinfo{year}{1928}) \bibinfo{pages}{32--74}.
	\bibitem[{Pierre et~al.(2019)Pierre, Nadia, Tan, Touzi, and
		Warin}]{HenryLabordere2019}
	\bibinfo{author}{H.~L. Pierre}, \bibinfo{author}{O.~Nadia},
	\bibinfo{author}{X.~Tan}, \bibinfo{author}{N.~Touzi},
	\bibinfo{author}{X.~Warin},
	\newblock \bibinfo{title}{Branching diffusion representation of semilinear
		{PDE}s and {M}onte {C}arlo approximation},
	\newblock \bibinfo{journal}{Ann. Inst. H. Poincar{\'e} Probab. Statist.}
	\bibinfo{volume}{55} (\bibinfo{year}{2019}) \bibinfo{pages}{184--210}.
	\bibitem[{E et~al.(2019)E, Martin, Arnulf, and Thomas}]{weinan2019multilevel}
	\bibinfo{author}{W.~E}, \bibinfo{author}{H.~Martin},
	\bibinfo{author}{J.~Arnulf}, \bibinfo{author}{K.~Thomas},
	\newblock \bibinfo{title}{On multilevel {P}icard numerical approximations for
		high-dimensional nonlinear parabolic partial differential equations and
		high-dimensional nonlinear backward stochastic differential equations},
	\newblock \bibinfo{journal}{J. Sci. Comput.} \bibinfo{volume}{79}
	(\bibinfo{year}{2019}) \bibinfo{pages}{1534--1571}.
	\bibitem[{Mainini(2012)}]{mainini2012description}
	\bibinfo{author}{E.~Mainini},
	\newblock \bibinfo{title}{A description of transport cost for signed measures},
	\newblock \bibinfo{journal}{J. Math. Sci.} \bibinfo{volume}{181}
	(\bibinfo{year}{2012}) \bibinfo{pages}{837--855}.
	\bibitem[{Ambrosio et~al.(2011)Ambrosio, Mainini, and
		Serfaty}]{ambrosio2011gradient}
	\bibinfo{author}{L.~Ambrosio}, \bibinfo{author}{E.~Mainini},
	\bibinfo{author}{S.~Serfaty},
	\newblock \bibinfo{title}{Gradient flow of the {Chapman-Rubinstein-Schatzman}
		model for signed vortices},
	\newblock \bibinfo{journal}{Ann. Inst. H. Poincar{\'e} Anal. Non lin{\'e}aire}
	\bibinfo{volume}{28} (\bibinfo{year}{2011}) \bibinfo{pages}{217--246}.
	\bibitem[{Piccoli et~al.(2019)Piccoli, Rossi, and
		Tournus}]{piccoli2019wasserstein}
	\bibinfo{author}{B.~Piccoli}, \bibinfo{author}{F.~Rossi},
	\bibinfo{author}{M.~Tournus},
	\newblock \bibinfo{title}{A {W}asserstein norm for signed measures, with
		application to nonlocal transport equation with source term},
	\newblock \bibinfo{journal}{arXiv preprint arXiv:1910.05105}
	(\bibinfo{year}{2019}).
	\bibitem[{Lawson(1967)}]{lawson1967generalized}
	\bibinfo{author}{J.~D. Lawson},
	\newblock \bibinfo{title}{Generalized {R}unge-{K}utta processes for stable
		systems with large {L}ipschitz constants},
	\newblock \bibinfo{journal}{SIAM J. Numer. Anal.} \bibinfo{volume}{4}
	(\bibinfo{year}{1967}) \bibinfo{pages}{372--380}.
	\bibitem[{Li et~al.(2016)Li, Yang, and Wong}]{LiYangWong2016}
	\bibinfo{author}{D.~Li}, \bibinfo{author}{K.~Yang}, \bibinfo{author}{W.~Wong},
	\newblock \bibinfo{title}{Density estimation via discrepancy based adaptive
		sequential partition},
	\newblock \bibinfo{journal}{Adv. Neural. Inf. Process. Syst.}
	(\bibinfo{year}{2016}) \bibinfo{pages}{1091--1099}.
	\bibitem[{E et~al.(2017)E, Han, and Jentzen}]{han2017deep}
	\bibinfo{author}{W.~E}, \bibinfo{author}{J.~Han}, \bibinfo{author}{A.~Jentzen},
	\newblock \bibinfo{title}{Deep learning-based numerical methods for
		high-dimensional parabolic partial differential equations and backward
		stochastic differential equations},
	\newblock \bibinfo{journal}{Commun. Math. Stat.} \bibinfo{volume}{5}
	(\bibinfo{year}{2017}) \bibinfo{pages}{349--380}.
	\bibitem[{N{\"u}sken and Richter(2023)}]{nusken2021interpolating}
	\bibinfo{author}{N.~N{\"u}sken}, \bibinfo{author}{L.~Richter},
	\newblock \bibinfo{title}{Interpolating between {BSDE}s and {PINN}s: deep
		learning for elliptic and parabolic boundary value problems},
	\newblock \bibinfo{journal}{J. Mach. Learn.} \bibinfo{volume}{2}
	(\bibinfo{year}{2023}) \bibinfo{pages}{31--64}.
	\bibitem[{Gao and Wang(2023)}]{GaoWang2023}
	\bibinfo{author}{W.~Gao}, \bibinfo{author}{C.~Wang},
	\newblock \bibinfo{title}{Active learning based sampling for high-dimensional
		nonlinear partial differential equations},
	\newblock \bibinfo{journal}{J. Comput. Phys.} \bibinfo{volume}{475}
	(\bibinfo{year}{2023}) \bibinfo{pages}{111848}.
	\bibitem[{Hochbruck and Ostermann(2010)}]{HochbruckOstermann2010}
	\bibinfo{author}{M.~Hochbruck}, \bibinfo{author}{A.~Ostermann},
	\newblock \bibinfo{title}{Exponential integrators},
	\newblock \bibinfo{journal}{Acta Numer.} \bibinfo{volume}{19}
	(\bibinfo{year}{2010}) \bibinfo{pages}{209--286}.
	\bibitem[{Dimov(2008)}]{bk:Dimov2008}
	\bibinfo{author}{I.~T. Dimov}, \bibinfo{title}{Monte Carlo Methods for Applied
		Scientists}, \bibinfo{publisher}{World Scientific}, \bibinfo{year}{2008}.
	\bibitem[{Kwa{\'s}nicki(2017)}]{kwasnicki2017ten}
	\bibinfo{author}{M.~Kwa{\'s}nicki},
	\newblock \bibinfo{title}{Ten equivalent definitions of the fractional
		{L}aplace operator},
	\newblock \bibinfo{journal}{Fract. Calc. Appl. Anal.} \bibinfo{volume}{20}
	(\bibinfo{year}{2017}) \bibinfo{pages}{7--51}.
	\bibitem[{Shao and Xiong(2020)}]{ShaoXiong2020}
	\bibinfo{author}{S.~Shao}, \bibinfo{author}{Y.~Xiong},
	\newblock \bibinfo{title}{Branching random walk solutions to the {W}igner
		equation},
	\newblock \bibinfo{journal}{SIAM J. Numer. Anal.} \bibinfo{volume}{58(5)}
	(\bibinfo{year}{2020}) \bibinfo{pages}{2589--2608}.
	\bibitem[{Raviart(1983)}]{raviart1983particle}
	\bibinfo{author}{P.~A. Raviart},
	\newblock \bibinfo{title}{An analysis of particle methods},
	\newblock in: \bibinfo{booktitle}{Numerical Methods in Fluid Dynamics. Vol.
		1127 of Lecture Notes in Mathematics.}, \bibinfo{publisher}{Springer},
	\bibinfo{year}{1983}, pp. \bibinfo{pages}{243--324}.
	\bibitem[{Yan and Caflisch(2015)}]{Yan2015}
	\bibinfo{author}{B.~Yan}, \bibinfo{author}{R.~E. Caflisch},
	\newblock \bibinfo{title}{{A Monte Carlo method with negative particles for
			{C}oulomb collisions}},
	\newblock \bibinfo{journal}{J. Comput. Phys.} \bibinfo{volume}{298}
	(\bibinfo{year}{2015}) \bibinfo{pages}{711--740}.
	\bibitem[{Wu et~al.(2017)Wu, Shin, and Xiu}]{wu2017randomized}
	\bibinfo{author}{K.~Wu}, \bibinfo{author}{Y.~Shin}, \bibinfo{author}{D.~Xiu},
	\newblock \bibinfo{title}{A randomized tensor quadrature method for high
		dimensional polynomial approximation},
	\newblock \bibinfo{journal}{SIAM J. Sci. Comput.} \bibinfo{volume}{39}
	(\bibinfo{year}{2017}) \bibinfo{pages}{A1811--A1833}.
	\bibitem[{Wu and Xiu(2018)}]{wu2018sequential}
	\bibinfo{author}{K.~Wu}, \bibinfo{author}{D.~Xiu},
	\newblock \bibinfo{title}{Sequential function approximation on arbitrarily
		distributed point sets},
	\newblock \bibinfo{journal}{J. Comput. Phys.} \bibinfo{volume}{354}
	(\bibinfo{year}{2018}) \bibinfo{pages}{370--386}.
	\bibitem[{Roberts and Rosenthal(2001)}]{Roberts2001metropolis}
	\bibinfo{author}{G.~O. Roberts}, \bibinfo{author}{J.~S. Rosenthal},
	\newblock \bibinfo{title}{Optimal scaling for various {M}etropolis-{H}astings
		algorithms},
	\newblock \bibinfo{journal}{Statistical Science} \bibinfo{volume}{16}
	(\bibinfo{year}{2001}) \bibinfo{pages}{351--367}.
	\bibitem[{Beylkin et~al.(1998)Beylkin, Keiser, and Vozovoi}]{beylkin1998new}
	\bibinfo{author}{G.~Beylkin}, \bibinfo{author}{J.~M. Keiser},
	\bibinfo{author}{L.~Vozovoi},
	\newblock \bibinfo{title}{A new class of time discretization schemes for the
		solution of nonlinear {PDE}s},
	\newblock \bibinfo{journal}{J. Comput. Phys.} \bibinfo{volume}{147}
	(\bibinfo{year}{1998}) \bibinfo{pages}{362--387}.
	\bibitem[{Hansen et~al.(2012)Hansen, Kramer, and Ostermann}]{hansen2012second}
	\bibinfo{author}{E.~Hansen}, \bibinfo{author}{F.~Kramer},
	\bibinfo{author}{A.~Ostermann},
	\newblock \bibinfo{title}{A second-order positivity preserving scheme for
		semilinear parabolic problems},
	\newblock \bibinfo{journal}{Appl. Numer. Math.} \bibinfo{volume}{62}
	(\bibinfo{year}{2012}) \bibinfo{pages}{1428--1435}.
	\bibitem[{Crouseilles et~al.(2020)Crouseilles, Einkemmer, and
		Massot}]{crouseilles2020exponential}
	\bibinfo{author}{N.~Crouseilles}, \bibinfo{author}{L.~Einkemmer},
	\bibinfo{author}{J.~Massot},
	\newblock \bibinfo{title}{Exponential methods for solving hyperbolic problems
		with application to collisionless kinetic equations},
	\newblock \bibinfo{journal}{J. Comput. Phys.} \bibinfo{volume}{420}
	(\bibinfo{year}{2020}) \bibinfo{pages}{109688}.
	\bibitem[{Gnewuch et~al.(2012)Gnewuch, Wahlström, and
		Winzen}]{GnewuchWahlstromWinzen2012}
	\bibinfo{author}{M.~Gnewuch}, \bibinfo{author}{M.~Wahlström},
	\bibinfo{author}{C.~Winzen},
	\newblock \bibinfo{title}{A new randomized algorithm to approximate the star
		discrepancy based on threshold accepting},
	\newblock \bibinfo{journal}{SIAM J. Numer. Anal.} \bibinfo{volume}{50(2)}
	(\bibinfo{year}{2012}) \bibinfo{pages}{781--807}.
	\bibitem[{Xiong and Shao(2019)}]{XiongShao2019}
	\bibinfo{author}{Y.~Xiong}, \bibinfo{author}{S.~Shao},
	\newblock \bibinfo{title}{The {W}igner branching random walk: {E}fficient
		implementation and performance evaluation},
	\newblock \bibinfo{journal}{Commun. Comput. Phys.} \bibinfo{volume}{25}
	(\bibinfo{year}{2019}) \bibinfo{pages}{871--910}.
	\bibitem[{Silverman(2018)}]{bk:Silverman2018}
	\bibinfo{author}{B.~W. Silverman}, \bibinfo{title}{{Density Estimation for
			Statistics and Data Analysis}}, \bibinfo{publisher}{Routledge},
	\bibinfo{year}{2018}.
	\bibitem[{Bungartz and Griebel(2004)}]{Griebel2004Sparsegrids}
	\bibinfo{author}{H.-J. Bungartz}, \bibinfo{author}{M.~Griebel},
	\newblock \bibinfo{title}{Sparse grids},
	\newblock \bibinfo{journal}{Acta Numer.} \bibinfo{volume}{13}
	(\bibinfo{year}{2004}) \bibinfo{pages}{147--269}.
	\bibitem[{Kormann and Sonnendr{\"u}cker(2016)}]{kormann2016sparse}
	\bibinfo{author}{K.~Kormann}, \bibinfo{author}{E.~Sonnendr{\"u}cker},
	\newblock \bibinfo{title}{Sparse grids for the {V}lasov-{P}oisson equation},
	\newblock in: \bibinfo{booktitle}{Sparse Grids and Applications-Stuttgart
		2014}, \bibinfo{publisher}{Springer}, \bibinfo{year}{2016}, pp.
	\bibinfo{pages}{163--190}.
	\bibitem[{Vershynin(2018)}]{Vershynin_2018}
	\bibinfo{author}{R.~Vershynin}, \bibinfo{title}{{H}igh-{D}imensional
		{P}robability: {A}n {I}ntroduction with {A}pplications in {D}ata {S}cience},
	\bibinfo{publisher}{Cambridge University Press}, \bibinfo{year}{2018}.
	\bibitem[{Ye et~al.(2024)Ye, Tian, and Wang}]{YeTianWang2024}
	\bibinfo{author}{Q.~Ye}, \bibinfo{author}{X.~Tian}, \bibinfo{author}{D.~Wang},
	\newblock \bibinfo{title}{A fast and accurate solver for the fractional
		{F}okker-{P}lanck equation with {D}irac-{D}elta initial conditions},
	\newblock \bibinfo{journal}{arXiv preprint arXiv:2407.15315}
	(\bibinfo{year}{2024}).
	
\end{thebibliography}

\appendix
\setcounter{equation}{0}
\setcounter{algorithm}{0}
\renewcommand\theequation{A.\arabic{equation}}
\renewcommand\thealgorithm{A.\arabic{algorithm}}
\section{A random walk algorithm for nonlocal Laplacian operator}
\label{appendix A}
Starting from the heat semigroup definition \eqref{fractional heat}, we can obtain the decomposition by introducing a small $\varepsilon$, 
\begin{equation}\label{eq:a1}
-(-\Delta)^{\alpha/2} g(\bx, t) = \frac{1}{|\Gamma(-\alpha/2)|} \left(\int_0^{\varepsilon} + \int_{\varepsilon}^{\infty}\right) (\me^{s \Delta} g(\bx, t) - g(\bx, t) ) \frac{\D s}{s^{1 + \alpha/2}}.
\end{equation}
For the first RHS term of Eq.~\eqref{eq:a1}, we can obtain an expansion utilizing $\me^x - 1 = x+\mathcal{O}(x^2)$,
\begin{equation}\label{first term}
	\frac{1}{|\Gamma(-\alpha/2)|} \int_0^{\varepsilon} \frac{\me^{s \Delta} - 1}{s^{1 + \alpha/2}} g(\bx, t) \D s = C(\alpha, \varepsilon) \Delta g(\bx,t) + \mathcal{O}(\varepsilon^{2 - \frac{\alpha}{2}}),
\end{equation}
and for the second RHS term of Eq.~\eqref{eq:a1}, a direct calculation gives
\begin{equation*}
	\frac{1}{|\Gamma(-\alpha/2)|} \int_{\varepsilon}^{\infty} (\me^{s \Delta} g(\bx, t) - g(\bx, t) ) \frac{\D s}{s^{1 + \alpha/2}} = \frac{1}{|\Gamma(-\alpha/2)|} \int_{\varepsilon}^{\infty}  \frac{\me^{s \Delta} g(\bx, t) }{s^{1 + \alpha/2}} \D s - \gamma(\alpha, \varepsilon)g(\bx,t),
\end{equation*}
where 
\begin{equation*}
	C(\alpha, \varepsilon) = \frac{\varepsilon^{1 - {\alpha}/{2}}}{(1 - {\alpha}/{2}) |\Gamma(-\alpha/2)|} > 0, \quad  \gamma(\alpha, \varepsilon) = \frac{\varepsilon^{-\alpha/2}}{(\alpha/2) |\Gamma(-\alpha/2)|} > 0.
\end{equation*}
After ignoring the higher-order term $\mathcal{O}(\varepsilon^{2 - \frac{\alpha}{2}})$ in Eq.~\eqref{first term}, we turn to solve a modified subdiffusion equation
\begin{equation}\label{modifiedl_subdiffusion}
	\frac{\partial }{\partial t}g(\bx, t) = C(\alpha, \varepsilon)  \Delta g(\bx, t) - \gamma(\alpha, \varepsilon) g(\bx, t) + \frac{1}{|\Gamma(-\alpha/2)|} \int_{\varepsilon}^{\infty}  \frac{\me^{s \Delta} g(\bx, t) }{s^{1 + \alpha/2}} \D s,\quad g(\bx, 0) = \varphi(\bx),
\end{equation}
and its integral formulation reads that  
\begin{equation}\label{g(x,tau)}
	g(\bx, \tau) =  \me^{-\gamma(\alpha, \varepsilon) \tau}\me^{C(\alpha, \varepsilon) \tau \Delta }\varphi(\bx) + \int_0^{\tau} \D t_1 ~ p(t_1) \int_{\varepsilon}^{\infty} \D s_1 ~ q(s_1) \, \me^{(C(\alpha, \varepsilon) t_1 + s_1) \Delta} g(\bx, \tau - t_1)
\end{equation}
with probability density functions
\begin{align}
	p(t) &= \gamma(\alpha, \varepsilon) \me^{-\gamma(\alpha, \varepsilon) t},\quad t\in [0, \infty), \label{pro p(t)} \\
	q(s) &= \frac{\alpha/2}{\varepsilon^{-\alpha/2}} \frac{1}{s^{1 + \alpha/2}},\quad s\in [\varepsilon, \infty). \label{pro q(s)}
\end{align}
In consequence, we can obtain a similar expression to Eq.~\eqref{g(x,tau)} for
$g(\bx, \tau - t_1)$, insert it back to Eq.~\eqref{g(x,tau)}, repeat this procedure and finally arrive at the Neumann series expansion,
\begin{equation}
	\begin{split}
		&g(\bx, \tau) = \me^{-\gamma(\alpha, \varepsilon) \tau + C(\alpha, \varepsilon) \tau \Delta}\varphi(\bx) \\
		&+ \int_0^{\tau} \D t_1~ p(t_1) \int_{\varepsilon}^{\infty} \D s_1 ~ q(s_1) \me^{(C(\alpha, \varepsilon) \tau + s_1) \Delta} \me^{-\gamma(\alpha,\epsilon)(\tau-t_1)}\varphi(\bx)\\
		&+ \int_0^{\tau} \D t_1 \,p(t_1)\int_0^{\tau-t_1} \D t_2 \,p(t_2)~  \int_{\varepsilon}^{\infty} \D s_1~ q(s_1) \int_{\varepsilon}^{\infty}   \D s_2~ q(s_2) \me^{(C(\alpha, \varepsilon) \tau + s_1 + s_2) \Delta}\me^{-\gamma(\alpha,\epsilon)(\tau-t_1-t_2)}\varphi(\bx)\\
		& + \cdots,
	\end{split}
\end{equation}
which yields a random walk algorithm for the nonlocal Laplacian operator in Example~\ref{example nonlocal} (see Algorithm~\ref{alg fractional random walk}). According to Bernstein's inequality \cite{Vershynin_2018} of the exponential density $p(t)$, we can obtain the probability that the while loop has not been exited at the $n$-th step in Algorithm~\ref{alg fractional random walk},
\begin{equation}\label{exp concen}
	 \mathbb{P}\left(\sum_{i=1}^n t_i \leq \tau\right) 
				\leq  2 \exp \left[-C_1 \min \left(\frac{(n-\tau \gamma(\alpha, \varepsilon))^2}{n C_2}, n-\tau \gamma(\alpha, \varepsilon)\right)\right],
\end{equation}
where $C_1, C_2>0$ are absolute constants. The probability in \eqref{exp concen} approaches zero exponentially with respect to $n$, and thus in practical simulations, the while loop in Algorithm~\ref{alg fractional random walk} terminates in a finite number of steps.

\begin{algorithm}[H]
	\caption{A random walk algorithm for the nonlocal operator $\mathcal{L} = -(-\Delta)^{\alpha / 2}$ in Example~\ref{example nonlocal}.} 
	\hspace*{0.02in} {\bf Input:}
	\label{alg fractional random walk}
	The power $\alpha$, the cut-off parameter $\epsilon$, the time step $\tau$ and the particle location $\bx_i(t_m)$. \\
	\hspace*{0.02in} {\bf Output:}  $\bx_i(t_{m+1})$ under the action of $\mathcal{L}$.
	\begin{algorithmic}[1]
		\State $n=1$;
		\While{true}
		\State Generate a random number $t_n$ that obeys the probability density $p(t)$ in Eq.~\eqref{pro p(t)};
		\If{$\sum_{i=1}^{n}t_i > \tau$}
		\State Break;
		\Else
		\State Generate a random number $s_n$ that obeys the probability density $q(s)$ in Eq.~\eqref{pro q(s)};
		\EndIf
		\State $n\leftarrow n+1$;
		\EndWhile
		\State $S = C(\alpha,\epsilon)\tau + \sum_{i=1}^{n-1}s_i$;
		\State Generate a random number $\by$ that obeys the normal distribution $\mathcal{N}(\vec{0}, 2S \vec{I})$;
		\State $\vec{x}_i(t_{m+1}) \leftarrow \vec{x}_i(t_m) - \vec{y}$;
	\end{algorithmic}
\end{algorithm}

Next we show the numerical experiments using SPM to integrate the high-dimensional linear nonlocal equation. Consider the nonlocal Fokker-Planck equation with Dirac-Delta initial conditions
\begin{equation}\label{FFPE}
	\begin{aligned}
		\frac{\partial}{\partial t} g(\bx, t) &= \bm{b} \cdot \nabla g(\bx, t) + c \Delta g(\bx, t) - (-\Delta)^{\alpha/2} g(\bx, t), \quad \bx\in \mathbb{R}^d, \\
		g(\bx, 0) &= \delta(\bx - \bx_0).
	\end{aligned}
\end{equation}
The algorithm for solving Eq.~\eqref{FFPE} using SPM is summarized in Algorithm~\ref{alg hd frac linear}. We can simulate particles one by one in Algorithm~\ref{alg hd frac linear} and the storage cost is only $\mathcal{O}(d)$. The fundamental solution of Eq.~\eqref{FFPE} can be represented as a one-dimensional integral in radial coordinate ($d\ge 2$) \cite{YeTianWang2024}
\begin{equation}\label{rad sol}
	g(\bx, t) = \widetilde g(y, t)= \frac{1}{y^{(d-2)/2}} \int_0^{\infty} \left(\frac{r}{2\pi}\right)^{d/2} J_{(d-2)/2}(yr) \exp\left(- (c r^2 + r^{\alpha}) t\right) \D r,
\end{equation}
where $y = |\bx - \bx_0 - \bm{b} t|$ and $J_{\nu}(x)$ is the Bessel function of the first kind. It is noted that the 2-D projection $M(x_1, x_2, t) = \int_{\mathbbm{R}^{d-2}} g(\bx,t) \D x_3 \dots \D x_d$ satisfies Eq.~\eqref{rad sol} for $d = 2$. The reference solution can be produced by a deterministic solver using the Laguerre-Gauss quadrature.

\begin{algorithm}[H]
	\caption{Stochastic particle method for nonlocal Fokker-Planck equation \eqref{FFPE}.} 
	\hspace*{0.02in} {\bf Input:}
	\label{alg hd frac linear}
	The cut-off parameter $\epsilon$, the final time $T$, the time step $\tau$, the particle number $N$ and the equation parameters $\vec{b}, c, \alpha, d, \bx_0$. \\
	\hspace*{0.02in} {\bf Output:}  $\langle \varphi(\bx), g(\bx, T)\rangle$ for any given test function $\varphi \in \mathbb{H}^2(\mathbb{R}^d)$.
	\begin{algorithmic}[1]
		\For{$i = 1:N$}
			\State $\bx_i(0) = \bx_0$, $w_i(0) = 1$;
			\For{$m = 1:\frac{T}{\tau}$}
			\State Simulate convection for gradient operator $\vec{b}\cdot \nabla$ according to Eq.~\eqref{simu convec};
			\State Simulate diffusion for Laplace operator $c\Delta$ according to Eq.~\eqref{simu diff};
			\State Simulate random walk for nonlocal operator $\mathcal{L} = -(-\Delta)^{\alpha / 2}$ according to Algorithm~\ref{alg fractional random walk};
			\EndFor
		\EndFor
		\State $\langle \varphi(\bx), g(\bx, T)$ = $\frac{1}{N} \sum_{i=1}^{N} w_i(T) \varphi(\bx_i(T))$;
	\end{algorithmic}
\end{algorithm}

We choose $T = 4$, $\tau = 0.1$, $\epsilon = 0.005$, $\alpha =1.5$, $\vec{b} = \left(1,1,\dots,1\right)$, $c=0.2$ and $d = 1000000$. In Algorithm~\ref{alg hd frac linear}, each particle is independent of each other and is inherently parallelizable. We employ 1000 cores for MPI parallelization. By taking the test function $\varphi(\bx) = x_1$ and $\varphi(\bx) = x_1^2$ respectively, the two macroscopic quantities 
\begin{equation}
	\begin{aligned}
		O_1(t) &= \int_{\mathbb{R}^{d}} x_1 u(\bx,t) \D \bx, \\
		O_2(t) &= \int_{\mathbb{R}^{d}} x_1^2 u(\bx,t) \D \bx,
	\end{aligned}
\end{equation}
are employed to measure the accuracy.

\begin{figure}[htbp]
	\centering
	\begin{subfigure}[b]{0.45\textwidth}
		\centering
		\includegraphics[width=\textwidth]{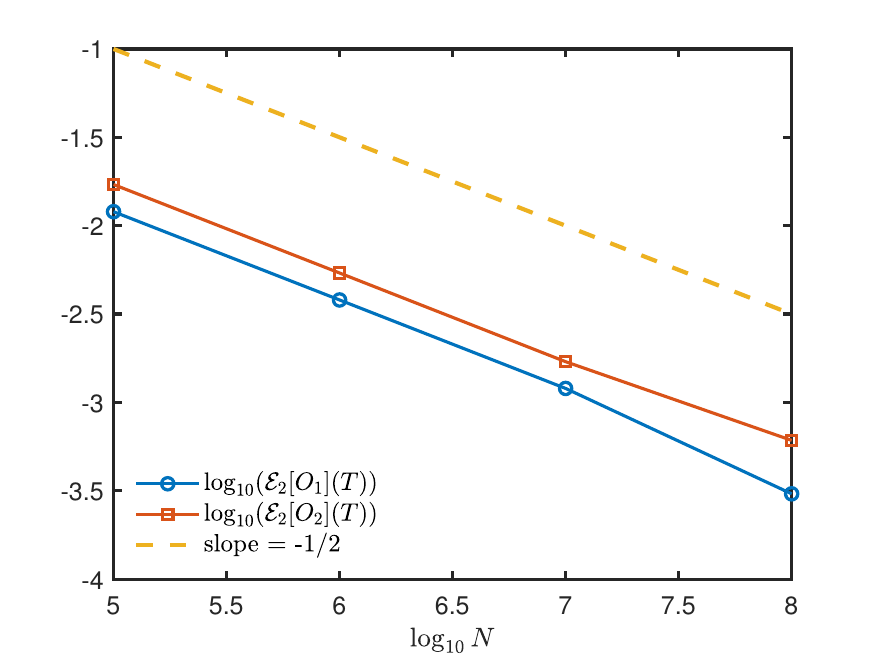}
	\end{subfigure}
	\caption{The 1000000-D linear nonlocal equation: SPM exhibits a half order convergence. Here $T=4$.}
	\label{fig frac conver}
\end{figure}

\begin{figure}[H]
	\centering
	\begin{subfigure}[b]{0.49\textwidth}
		\centering
		\includegraphics[width=\textwidth]{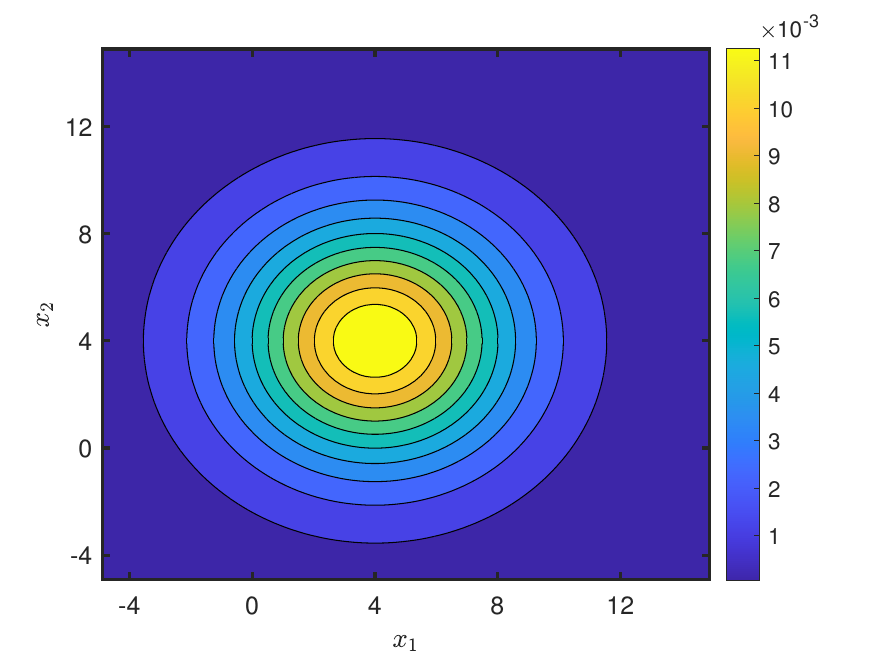}
		\caption{Reference solution.}
	\end{subfigure}
	\hfill
	\begin{subfigure}[b]{0.49\textwidth}
		\centering
		\includegraphics[width=\textwidth]{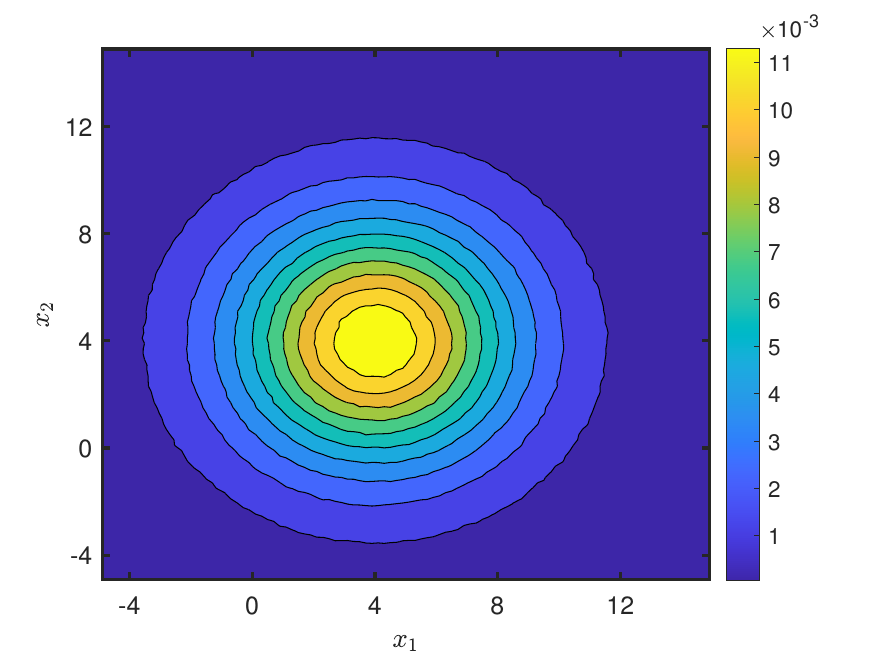}
		\caption{SPM solution with $N=1\times 10^8$.}
	\end{subfigure}
	\caption{The 1000000-D linear nonlocal equation: Filled contour plots of $M(x_1, x_2, 4)$ for (a) the reference solution and (b) the numerical solution produced by SPM with $N=1\times 10^8$.}
	\label{hdlinearfracfig}
\end{figure}

Table~\ref{FRW_convergence} provides the relative errors $\mathcal{E}_2[O_1](T)$ and $\mathcal{E}_2[O_2](T)$ under different sample sizes $N$. The convergence is verified in Figure~\ref{fig frac conver}. Figure~\ref{hdlinearfracfig} plots 2-D projections $M(x_1, x_2, T)$ for the numerical solution produced by SPM against the reference solution. A good agreement is observed.

\begin{table}[htbp]
	\centering
	\begin{tabular}{cccc}
		\toprule
		 $N$ & $\mathcal{E}_2[O_1](T)$ & $\mathcal{E}_2[O_2](T)$  & Time/s\\
		\midrule
	
			$1\times10^5$	&	0.0120		&	0.0171		&	7.81 	\\
				$1\times10^6$	&	0.0038		&	0.0054		&	77.79 	\\
				$1\times10^7$	&	0.0012		&	0.0017		&	771.01	\\ 
				$1\times10^8$	&	0.0003		&	0.0006		&	7970.63	\\
				
		\bottomrule
	\end{tabular}  
	\caption{The 1000000-D linear nonlocal equation: Total wall time with 1000 cores in parallel and the relative errors under different sample sizes at $T=4$.}
	\label{FRW_convergence}
\end{table}

\end{document}